\documentclass[a4paper,12pt,reqno]{amsart}

\usepackage{amsmath, amsfonts, amssymb, amsthm, mathrsfs, mathtools, mathalfa}
\usepackage{comment}
\usepackage{yfonts}
\usepackage{enumerate}
\usepackage{enumitem}
\usepackage{stmaryrd}
\usepackage{fancyhdr}
\usepackage{bm}
\usepackage[utf8]{inputenc}
\usepackage[pdfencoding=auto, hyperfootnotes=false]{hyperref}
\usepackage{tikz,tikz-cd}
\usepackage[hang,flushmargin]{footmisc}
\usepackage{graphicx}
\usepackage{xr}
\usepackage{extarrows}
\usepackage[utf8]{inputenc}
\usepackage{subfig}

\newtheorem{theorem}{Theorem}

\newtheorem{property}[theorem]{Property}

\begin{document}

\title{The edge labeling of higher order Voronoi diagrams}

\author[M. Claverol et al.]{Mercè Claverol}
\address{Departament de Matemàtiques, Universitat Politècnica de Catalunya, Barcelona, Spain.}
\email{merce.claverol@upc.edus}

\author[]{Andrea de las Heras Parrilla}
\address{Departament de Matemàtiques, Universitat Politècnica de Catalunya, Barcelona, Spain.}
\email{andrea.de.las.heras@estudiantat.upc.edu}

\author[]{Clemens Huemer}
\address{Departament de Matemàtiques, Universitat Politècnica de Catalunya, Barcelona, Spain.}
\email{clemens.huemer@upc.edu}

\author[]{Alejandra Martínez-Moraian}
\address{Departamento de Física y Matemáticas\\ Universidad de Alcalá \\ Madrid\\ Spain}
\email{alejandra.martinezm@uah.es}


\begin{abstract}
We present an edge labeling of order-$k$ Voronoi diagrams, $V_k(S)$, of point sets $S$ in the plane, and study properties of the regions defined by them. Among them, we show that $V_k(S)$ has a small orientable cycle and path double cover, and we identify configurations that cannot appear in $V_k(S)$ for small values of $k$.
This paper also contains a systematic study of well-known and new properties of $V_k(S)$, all whose proofs only rely on elementary geometric arguments in the plane. The maybe most comprehensive study of structural properties of $V_k(S)$ was done by D.T. Lee (On k-nearest neighbor Voronoi diagrams in the plane) in 1982. Our work reviews and extends the list of properties of higher order Voronoi diagrams.
\end{abstract}

\maketitle

\section{Introduction}\label{sec1}
\begin{figure}[h!]
 	\subfloat[$V_1(S)$]{
 	\noindent
 		\begin{minipage}[c][1\width]{
 				0.32\textwidth}
 			\centering
 			\includegraphics[scale=0.22,page=1]{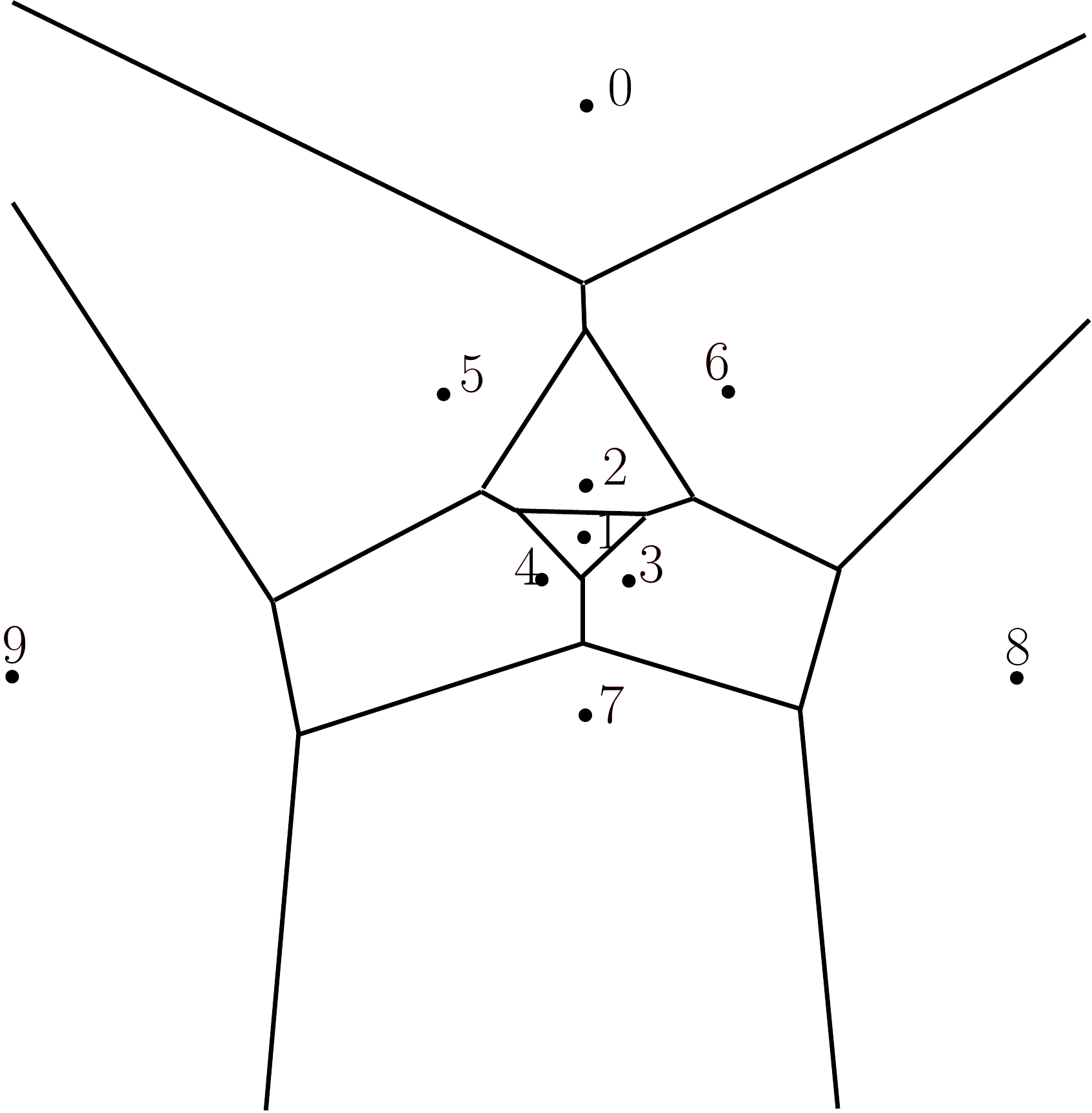}
 	\end{minipage}}
 	\hfill 	
 	\subfloat[$V_2(S)$]{
 	\noindent
 		\begin{minipage}[c][1\width]{
 				0.32\textwidth}
 			\centering
 			\includegraphics[scale=0.22,page=2]{FV1-V9.pdf}
 	\end{minipage}}
 	\hfill	
 	\subfloat[$V3(S)$]{
 	\noindent
 		\begin{minipage}[c][1\width]{
 				0.32\textwidth}
 			\centering
 			\includegraphics[scale=0.22,page=3]{FV1-V9.pdf}
 	\end{minipage}}
 	\hfill	
 	\subfloat[$V_4(S)$]{
 	\noindent
 		\begin{minipage}[c][1\width]{
 				0.32\textwidth}
 			\centering
 			\includegraphics[scale=0.22,page=4]{FV1-V9.pdf}
 	\end{minipage}}
 	\hfill	
 	\subfloat[$V_5(S)$]{
 	\noindent
 		\begin{minipage}[c][1\width]{
 				0.32\textwidth}
 			\centering
 			\includegraphics[scale=0.22,page=5]{FV1-V9.pdf}
 	\end{minipage}}
 	\hfill	
 	\subfloat[$V_6(S)$]{
 	\noindent
 	\noindent
 		\begin{minipage}[c][1\width]{
 				0.32\textwidth}
 			\centering
 			\includegraphics[scale=0.22,page=6]{FV1-V9.pdf}
 	\end{minipage}}
 	\hfill	
 	\subfloat[$V_7(S)$]{
 	\noindent
 		\begin{minipage}[c][1\width]{
 				0.32\textwidth}
 			\centering
 			\includegraphics[scale=0.22,page=7]{FV1-V9.pdf}
 	\end{minipage}}
 	\hfill	
 	\subfloat[$V_8(S)$]{
 	\noindent
 		\begin{minipage}[c][1\width]{
 				0.32\textwidth}
 			\centering
 			\includegraphics[scale=0.22,page=8]{FV1-V9.pdf}
 	\end{minipage}}
 	\hfill	
 	\subfloat[$V_9(S)$]{
 	\noindent
 		\begin{minipage}[c][1\width]{
 				0.32\textwidth}
 			\centering
 			\includegraphics[scale=0.22 ,page=9]{FV1-V9.pdf}
 	\end{minipage}}
 	\caption{Illustrating the central part of all Voronoi diagrams $V_k(S)$ for a set of points $S=\{0,\cdots, 9\}$. In each face, the $k$ nearest neighbors are indicated.}
 	\label{9VDs}
 \end{figure}

Let $S$ be a set of $n$ points in the plane such that no three of them are collinear and no four of them are cocircular, and let $1\leq k \leq n-1$ be an integer. The order-$k$ Voronoi diagram of $S$, $V_k(S)$, is a subdivision of the plane into faces such that the points in the same face have the same $k$ nearest points of $S$, also called $k$ nearest neighbors. See Figure~\ref{9VDs}.

Voronoi diagrams have applications in a broad range of disciplines, see e.g.~\cite{A91}. They are also known as Dirichlet tesselations or as Thiessen polygons~\cite{A91}.  The most studied Voronoi diagrams of point sets are $V_1(S)$, the classic Voronoi diagram, and $V_{n-1}(S)$,  the furthest point Voronoi diagram, which only has unbounded faces. Many properties of $V_k(S)$ were obtained by Lee~\cite{L82}, we also mention~\cite{DE84, E87,EI18,ESS11,L03,MRT19,OBSC2000,S98} among the sources on the structure of $V_k(S)$. In this work we review several of these structural results with new proofs, and we also present new results on $V_k(S).$ Most of them are based on a labeling of the edges of $V_k(S)$. An edge that delimits a face of $V_k(S)$ is a (possibly unbounded) segment of the perpendicular bisector of two points $i$ and $j$ of $S$. 
This well-known observation induces a natural labeling of the edges of $V_k(S)$ with the following rule:\\

\begin{figure}
	\centering
		\includegraphics[scale=0.7]{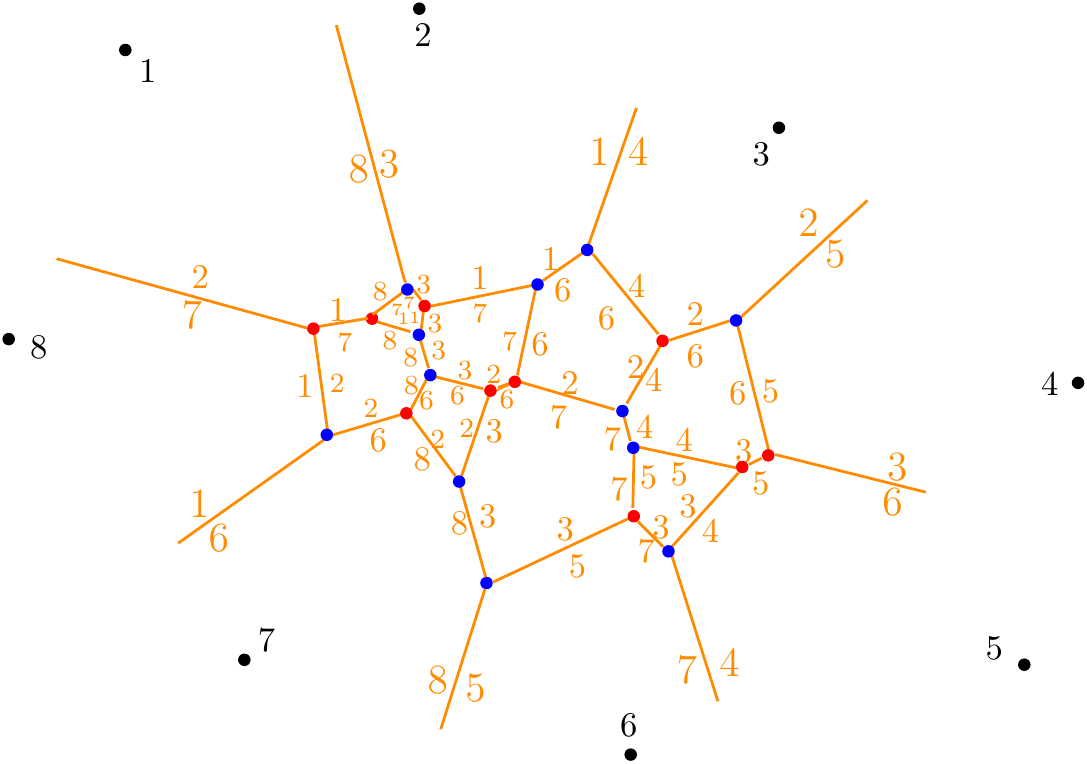}
	\caption{The edge labeling of $V_3(S)$ for a set $S$ of eight points in convex position. Vertices of type I are drawn in blue, and vertices of type II in red.}
	\label{fig:F-bisectors}
\end{figure}

$\bullet$ {\bf{Edge rule:}}
An edge of $V_k(S)$ which belongs to the perpendicular bisector of points $i,j \in S$ has labels $i$ and $j$, where we put the label $i$ on the side (half-plane) of the edge that contains point $i$ and we put label $j$ on the other side. See Figure~\ref{fig:F-bisectors}.\\

 Based on the edge rule, we also deduce a vertex rule and a face rule, proved in Section~\ref{sec:main}. 
$V_k(S)$ contains two types of vertices, denoted as type I and type II (also called new and old vertices~\cite{L82}), that are defined in Subsection \ref{Notation} of notation.  \\

$\bullet$ {\bf{Vertex rule:}}
Let $v$ be a vertex of $V_k(S)$ and let $\{i, j, \ell\}\in S$ be the set of labels of the edges incident to $v$. 
The cyclic order of the labels of the edges around $v$ is $i, i, j, j, \ell, \ell$ if $v$ is of type I, and it is $i, j, \ell, i, j, \ell$ if $v$ is of type II.\\

$\bullet$ {\bf{Face rule:}}
In each face of $V_k(S)$, the edges that have the same label $i$ are consecutive, and these labels $i$ are either all in the interior of the face, or are all in the exterior of the face.\\

Note that when walking along the boundary of a face, in its interior (exterior), a change in the labels
of its edges appears whenever we reach a vertex of type II (type I), see Figure~\ref{fig:F-bisectors}.

\begin{figure}[t]
    \centering
    \includegraphics[scale=0.7]{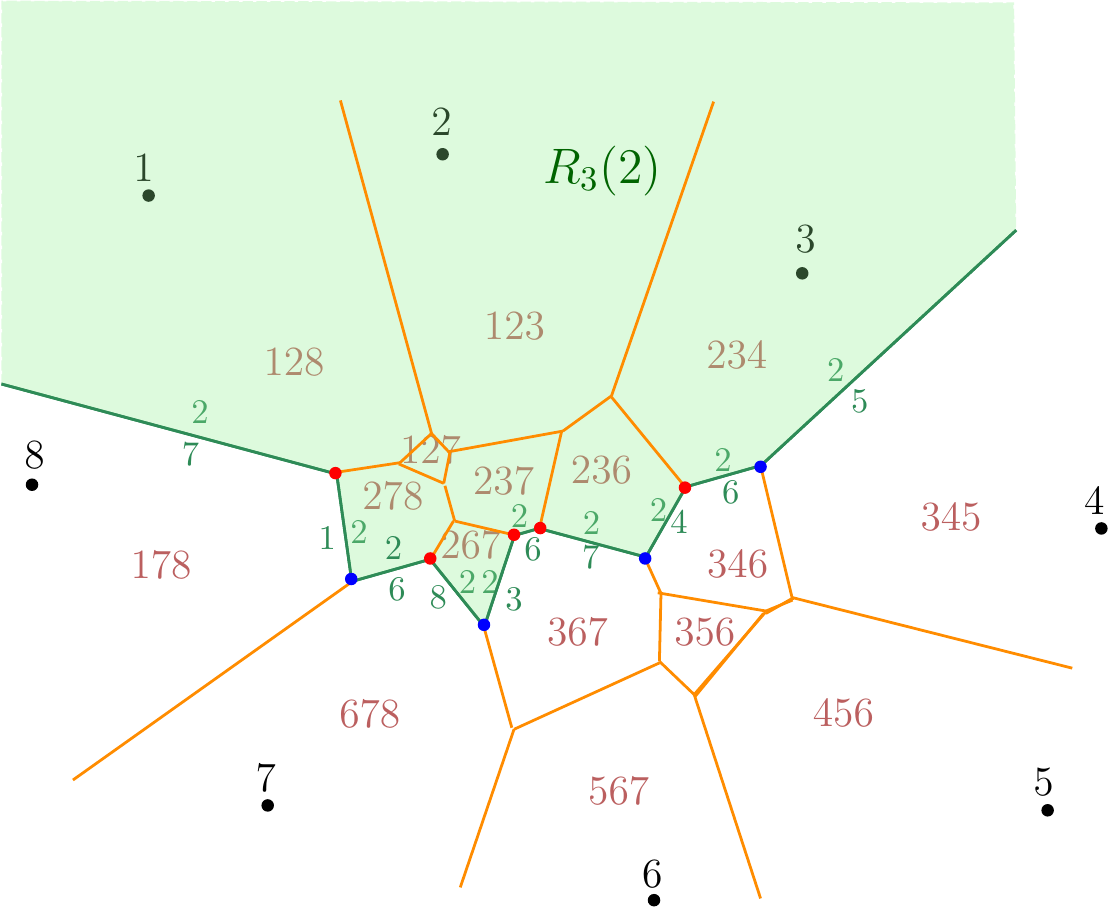}
    \caption{$V_3(S)$ for $S$ the point set in Figure~\ref{fig:F-bisectors}; in each face, its three nearest neighbors are indicated. In green, the region $R_3(2)$ formed by all the faces of $V_3(S)$ that have point $2$ as one of their three nearest neighbors. 
		The boundary of $R_3(2)$ is formed by all the edges which have the label $2$ and this label is always inside $R_3(2)$. The boundary vertices of $R_3(2)$ with an incident edge lying in the interior of $R_3(2)$ are of type II in $V_3(S)$ and the remaining boundary vertices are of type I in $V_3(S)$.}
    \label{fig:R32}
\end{figure}

We show that the edges with same label $i$ in $V_k(S)$ either form a cycle or a set of paths whose first and last edge are unbounded edges of $V_k(S)$. 
Edges with same label $i$ enclose a region $R_k(i)$ that consists of all the points of the plane that have point $i \in S$ as one of their $k$ nearest neighbors from $S$. See Figure~\ref{fig:R32}.
The union of all these regions $R_k(i)$ is a $k$-fold covering of the plane.  $R_k(i)$ is related to the $k$-th nearest point Voronoi diagram of $S$, that assigns to each point of the plane its $k$-nearest neighbor from $S$~\cite{OBSC2000}; this diagram is also called $k$-th degree
Voronoi diagram in~\cite{E87}. The region of a point $i \in S$ in the $k$-th nearest point Voronoi diagram is $R_k(i)\backslash R_{k-1}(i).$
A region of the $k$-th nearest point Voronoi diagram is not necessarily convex and may be disconnected~\cite{OBSC2000}. Such regions are also known as Brillouin zones, and their properties have been studied mainly for lattices, see e.g.~\cite{J84}, but also for arbitrary discrete sets~\cite{V00}. 

Edelsbrunner and Iglesias-Ham~\cite{EI18} proved that $R_k(i)$ is star-shaped. We further observe that $R_1(i)$ is contained in the kernel of this star-shaped set $R_k(i)$, and we identify the reflex  (convex) vertices on the boundary $B_k(i)$ of $R_k(i)$ as vertices of type II (type I).

We also show that every higher-order Voronoi diagram admits an {\it{orientable double cover}}~\cite{J85} of its edges using, precisely, the cycles and paths in $\cup_{i \in S} B_k(i).$
 A cycle and path double cover of a graph $G$ is a collection of cycles and paths $\mathcal{C}$ such that every edge of $G$ belongs to exactly two elements of $\mathcal{C}$. Paths are needed in a double cover $\mathcal{C}$ of $V_k(S)$ due to the unbounded edges. The number of paths in the obtained double cover $\mathcal{C}$ equals the number of unbounded edges of $V_k(S)$.
A double cover $\mathcal{C}$ is orientable if an orientation can be assigned
to each element of $\mathcal{C}$ such that for every edge $e$ of $G$,
the two cycles, resp. paths, that cover $e$ are oriented
in opposite directions through $e$~\cite{J85}.
It is well known that simple bridgeless planar graphs have an orientable cycle double cover whose cycles are formed by the edges bounding a face~\cite{B90}. Note that the double cover $\mathcal{C}$ of $V_1(S)$ formed by the boundaries $B_1(i)$ of $R_1(i)$, for $i \in S$, yields essentially the same double cover, the only difference being that unbounded regions of $V_1(S)$ are delimited by paths instead of cycles. 
The {\it{small  cycle double cover conjecture}} states that every simple bridgeless graph on $n$ vertices has a  cycle double cover with at most $n-1$ cycles, see Bondy~\cite{B90}. Seyffarth~\cite{S93} proved this conjecture for simple $4$-connected planar graphs, and also proved that any simple bridgeless planar graph of size $n$ has a cycle double cover with at most $3\lfloor(n-1)/2\rfloor$ cycles~\cite{S92}.
We show 
 that a higher order Voronoi diagram admits a much smaller orientable cycle and path double cover compared to its number of vertices.
For the cases when $k=\lfloor{n/2}\rfloor$ and $k=\lceil{n/2}\rceil$, we also show that $V_k(S)$ has a small orientable path double cover.
In addition, for point sets $S$ in convex position, we prove that $V_k(S)$ has an orientable path double cover consisting of $n$ paths, for any value of~$k$.

\begin{figure}[h!]
	\centering
	\subfloat[]{
		\includegraphics[scale=0.55,page=1]{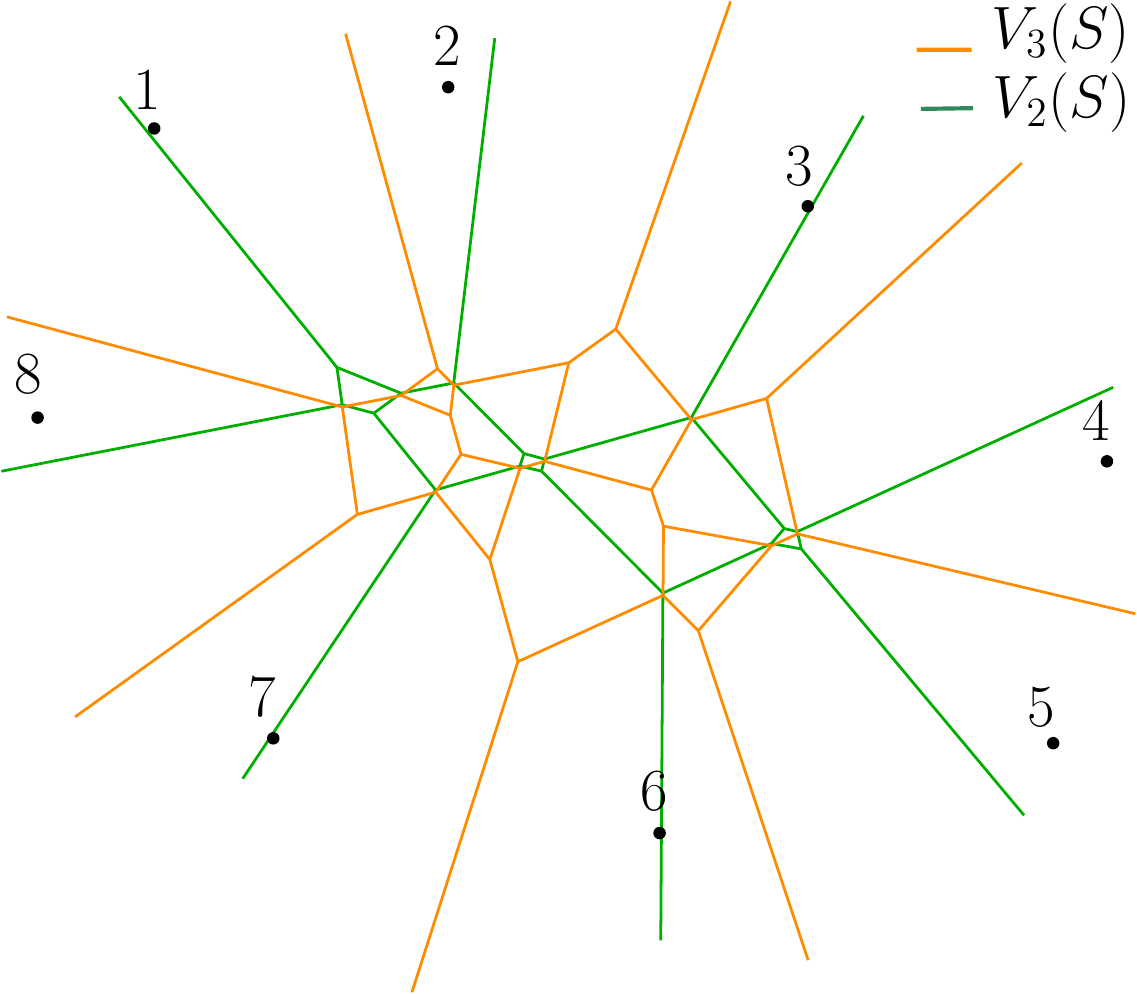}
		\label{fig:Th1}
	}~
	\subfloat[]{
		\includegraphics[scale=0.55,page=2]{FV23-34.pdf}
		\label{fig:Th2}
	}
	\caption{Consecutive higher order Voronoi diagrams for  the set $S=\{1,\cdots,8\}$ from Figure~\ref{fig:F-bisectors}.
		(a) $V_2(S)$ and $V_3(S)$.
		(b) $V_3(S)$ and $V_4(S)$.
	}
	\label{fig:FV234}
\end{figure}

We show several more new properties of $V_k(S)$, some dedicated to its unbounded faces. Also new proofs of known properties are given. For instance, it is known that the subgraph of $V_{k-1}(S)$ (also of $V_{k+1}(S)$\,) that lies inside a bounded face of $V_k(S)$ has the structure of a tree; see Figure~\ref{fig:FV234}. We reprove this fact using the edge labeling; see~\cite{E87,L82} for different proofs. As for unbounded faces, this only holds for $V_{k-1}(S)$; we show that the subgraph of $V_{k+1}(S)$ inside an unbounded face of $V_k(S)$ is not always a tree, but a forest. 

We also present a new proof for the number of vertices of $V_k(S)$, obtained by point moves; see~\cite{E87,L82,S98} for other proofs. Also tight upper bounds for the numbers of vertices of type I and of type II in bounded faces of $V_k(S)$ are shown; we could not find this result in the literature.

Finally, we use the edge labeling to show that certain configurations cannot appear in $V_k(S)$, for small values of $k$.\\

This paper is organized as follows. We start defining some notation in the following Subsection~\ref{Notation}. In Section~\ref{sec:main} we present a systematic study of basic properties of $V_k(S)$. The following sections rely on these basic properties.
In Section~\ref{sec:unbounded} we extend some of the properties for bounded faces to unbounded faces of $V_k(S)$. 
Section~\ref{sec:R_k} is on properties of the regions $R_k(i)$.
 In Section~\ref{sec:double} we present the double cover of $V_k(S)$.  Section~\ref{sec:count_V_k(S)} contains a proof for the number of vertices of $V_k(S)$. Finally, Section~\ref{sec:hexagons} is on alternating hexagons and on configurations that cannot appear in $V_k(S)$, for small values of $k$.\\

We give a proof for each property stated in this paper. Then, this work is self-contained, but we point out that many of the presented properties are well known. While some properties for $V_k(S)$ were proven by lifting points to $\mathbb{R}^d$, for $d>2$, and considering hyperplane arrangements, see e.g.~\cite{A90, E87, ES86}, all our proofs rely on elementary geometric arguments in the plane. We mention references for known alternative proofs in the text, and hope not to omit relevant references. Our exposition is purely based on structural properties of higher order Voronoi diagrams of point sets in the plane; we refer the reader to~\cite{LPL15} for references on algorithmic aspects, which are not treated in our text, and to~\cite{OBSC2000,A91} for surveys on the many variants of Voronoi diagrams.\\

\subsection{Notation}\label{Notation}

$S$ is a set of $n$ points in general position, that is, no three points of $S$ are collinear and no four points of $S$ are cocircular.
We denote the points of $S$ either by $\{p_1,\ldots,p_n\}$ or by $\{1,\ldots,n\}$.
The Voronoi diagram of order $k$ of $S$, $V_k(S)$, is a subdivision of the plane into faces, also called cells, such that the points in the same cell have the same $k$ nearest points of $S$, also called $k$ nearest neighbors, for $1 \leq k \leq n-1$. All distances refer to the Euclidean metric.  For points on the boundary of a cell, the $k$-th nearest neighbor of $S$ is not unique. In the literature, these cells are sometimes defined as open and sometimes as closed. In this text, faces of $V_k(S)$ are always closed.  We denote by $P_k$ a subset of $k$ points of $S$ defining a face $f(P_k)$ of $V_k(S)$. The points of $P_k$ are the neighbors in $V_k(S)$ of the points of $f(P_k)$. Note that not every subset of $k$ points of $S$ defines a face of $V_k(S)$. 

For each pair of points $i$, $j$ of $S$, let $\overline{ij}$ be the segment that connects $i$ and $j$, and let $b_{ij}$ be the perpendicular bisector of $\overline{ij}$. Each bisector $b_{ij}$ bounds two half-planes, $h(i,j)$ and $h(j,i)$, determining the points of the plane that are closer to $i$ than to $j$, and vice versa. Therefore, each face $f(P_k)$ of $V_k(S)$
is the following intersection of half-planes: $f(P_k)=\bigcap_{i\in P_k, j\in S\setminus P_k}h(i,j)$.

 Faces can be convex polygons or unbounded convex regions delimited by a polygonal chain. An unbounded face is denoted by $f^\infty(P_k)$ and a bounded face by $f^b(P_k)$. An edge on the boundary of two faces $f(P_k)$ and $f(Q_k)$ is a (possibly unbounded) segment on a bisector $b_{ij}$; it is denoted by $\overline{b_{ij}}$. Note that this is an abuse of notation, for several edges might be called $ \overline{b_{ij}}$. A vertex of $V_{k}(S)$ is the intersection point of three bisectors $b_{ab},b_{ac},$ and $b_{bc}$. Equivalently, it is the circumcenter of a circle $C_{abc}$ passing through the points $a,b,c$ of $S$. Such a vertex is denoted by $abc$. We distinguish two types of vertices in $f(P_k)$. If $a\in P_k$ and $b,c\in S\setminus P_k$ we say that $abc$ is of type I in $f(P_k)$; and if $a,b\in P_k$ and $c\in S\setminus P_k$, we say that it is of type II in $f(P_k)$. A vertex has the same type for its three incident faces, and hence a vertex of type I (type II) in $f(P_k)$ is a vertex of type I (type II) in $V_k(S)$. Note that a vertex $abc$ of $f(P_k)$ cannot be the center of a circle through three points of $P_k$. Lee~\cite{L82} and Lindenbergh~\cite{L03} refer to these vertices as {\it{new}} and {\it{old}} vertices, respectively. The latter notation is motivated by the property that each vertex of $V_k(S)$ is also a vertex either of $V_{k+1}(S)$ (new) or of $V_{k-1}(S)$ (old), see Property~\ref{property:alternate}. Instead, type I and type II indicates the number of points of $P_k$ in a vertex $abc$.
The closed disk bounded by $C_{abc}$ is denoted by $D_{abc}.$ With $c_k$  we denote the number of circles through three points of $S$ that enclose exactly $k$ other points of $S$.

We sometimes refer to $V_k(S)$ as a (drawing of a) graph, although it contains unbounded edges.
When considering the union of $V_k(S)$ and $V_{k+1}(S)$, the graph induced by $V_{k+1}(S)$ in a cell $f(P_k)$ of $V_k(S)$ is the subgraph of $V_{k+1}(S)$ whose vertices and edges are contained in $f(P_k)$. Similarly, the graph induced by $V_{k-1}(S)$ in $f(P_k)$ is the subgraph of $V_{k-1}(S)$ whose vertices and edges are contained in $f(P_k)$. 

The set of all the edges of $V_k(S)$ with label $i$ is denoted by $B_k(i)$. The region of the plane enclosed by $B_k(i)$ and that contains point $i$ is denoted by $R_k(i)$. Equivalently, $R_k(i)$ is the union of all the cells of $V_k(S)$ that have the point $i$ as one of their $k$ nearest neighbors.

For $0 \leq j \leq n-2$, a segment $\overline{pq}$ connecting two points $p$, $q$ of $S$ is a {\it{$j$-edge}} of $S$ if the oriented line from $p$ to $q$ divides the plane into two open half-planes, such that the half-plane to its left contains $j$ points of $S$. The number of $j$-edges of $S$ is denoted by $e_j$. 

\section{Properties of $V_k(S)$}\label{sec:main}

We start with basic properties of $V_k(S)$. Especially Properties \ref{property:adjacentvoronoi}-\ref{property:ija-ijb} are very elementary from the definition of $V_k(S)$. We state them explicitly, as they constitute the basis of our labeling rule, and they are used throughout the paper. Also most of the subsequent properties are well-known, see Lee~\cite{L82} and Edelsbrunner~\cite{E87}. In particular, see~\cite{L82} for Properties \ref{property:adjacentvoronoi}-\ref{property:alternate}. The statements of Properties \ref{property:connected}-\ref{property:nocycle} and \ref{prop:trees}-\ref{prop:treesize} referring to $V_{k-1}(S)$ also appear in \cite{L82}. Book \cite{E87} contains these statements for $V_{k+1}(S)$. We show that the analogous statement to Property~\ref{property:nocycle} for $V_{k+1}(S)$ is not true for unbounded faces; see Section~\ref{sec:unbounded}.  
Many proofs in~\cite{L82} rely on the construction of $V_k(S)$ from $V_{k-1}(S)$, whereas~\cite{E87} makes use of hyperplane arrangements. Some of the properties in this section have also been proved for abstract Voronoi diagrams, see~\cite{B15, B19}.

\begin{property}\label{property:adjacentvoronoi}
		Let $f(P_k)$ and $f(Q_k)$ be two faces of $V_k(S)$ sharing an edge. Then, $Q_k= \left(P_{k} \setminus \{i\}\right) \cup \{j\}$ for some $i, j\in S$. Furthermore, the edge shared by $f(P_k)$ and $f(Q_k)$ is a segment $\overline{b_{ij}}$ on the bisector $b_{ij}$. 
\end{property}
\begin{proof}
By definition of $V_k(S)$, there is at least one point $i\in P_k \setminus Q_k$ and one point $j\in Q_k\setminus P_k$. Then the common boundary of $f(P_k)$ and $f(Q_k)$ contains a segment of the bisector $b_{ij}$. Since faces of $V_k(S)$ are the intersection of half-planes, they are convex, and then the common boundary of $f(P_k)$ and $f(Q_k)$ can only contain one such segment. Hence $Q_k = \left(P_{k} \setminus \{i\}\right) \cup \{j\}$. See Figure~\ref{fig:R32}. 

 \end{proof}

\begin{property}\label{property:cell-edge}
	Let $\overline{b_{ij}}$ be an edge delimiting a face $f(P_k)$. If $i \in P_k$, then $i$ belongs to the half-plane defined by $b_{ij}$ that contains $f(P_k)$ and $j$ does not.
\end{property}
\begin{proof}
	By Property~\ref{property:adjacentvoronoi}, exactly one of $i$ or $j$ belongs to $P_k$. Since $b_{ij}$ is the perpendicular bisector of the points $i$ and $j$, exactly one of $i$ or $j$ lies in the half-plane defined by $b_{ij}$ that contains $f(P_k)$.
	Let $x$ be a point on  $b_{ij}$, then $x$ is equidistant to points $i$ and $j$. Moving $x$ to the interior of $f(P_k)$, $i$ gets closer to $x$ than $j$, by definition of $P_k$. Therefore, the point in the half-plane defined by $b_{ij}$ that contains $f(P_k)$ must be $i$.
	
\end{proof}

\begin{property}\label{property:ija-ijb}
	Any two consecutive vertices on a face $f(P_k)$ of $V_k(S)$ are of the form
	$ija$ and $ijb$.
	The edge connecting them is $\overline{b_{ij}}$.
\end{property}

\begin{figure}[h!]	
	\begin{center}
		\includegraphics[width=0.6\textwidth]{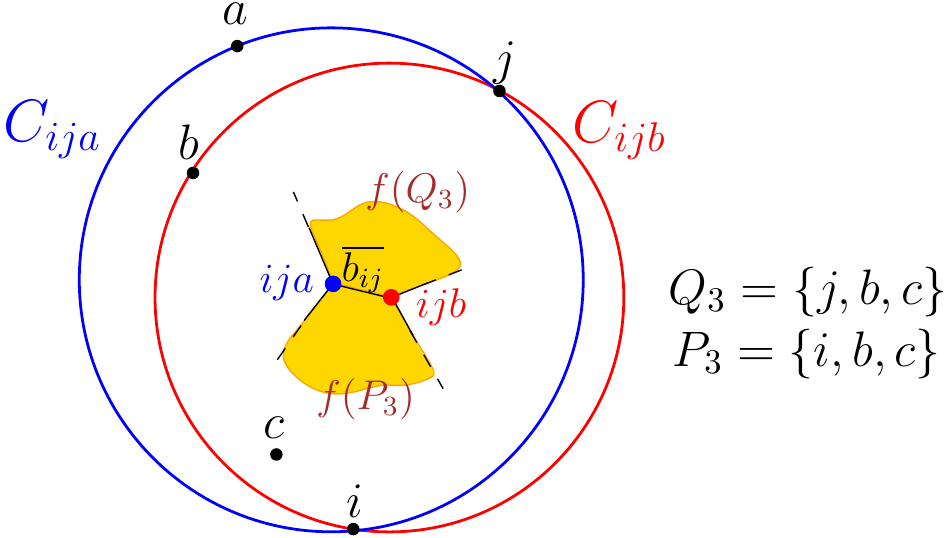}
		\caption{The edge $\overline{b_{ij}}$ connects vertices $ija$ and $ijb$. The sets $P_3$ and $Q_3$ differ in points $i,j$. The circle $C_{ijb}$ contains one point of $P_3$ (that is, $k-2$ points for $k=3$); then the vertex $ijb$ is of type II. The circle  $C_{ija}$ contains two (that is, $k-1$) points of $P_3$; then the vertex $ija$ is of type I.}
		\label{fig:Fobs1234}
	\end{center}
\end{figure}

\begin{proof}
	A vertex of $V_k(S)$ is the intersection point of three bisectors, which is the center of a circle passing through three points of $S$. Adjacent vertices are incident to a common bisector. See Figure~\ref{fig:Fobs1234}.
\end{proof}

\begin{property}\label{property:discs}
	Let $\overline{b_{ij}}$ be an edge delimiting a face $f(P_k)$. Then, $P_k \subset D_{ija} \cap D_{ijb}$, where $ija$ and $ijb$ are the endpoints of $\overline{b_{ij}}$.
\end{property}
\begin{proof}
  	The disk $D_{ija}$ centered in $ija$ and passing through the points $i$, $j$ and $a$ must contain the points of $P_k$ because $ija$ is a point of $f(P_k)$ and at least one of the points in $\left\lbrace i, j,a \right\rbrace$ is not in $P_k$. The same argument holds for $D_{ijb}$.
\end{proof}

\begin{property}\label{property:points-types-edges}
	The vertices of a face $f(P_k)$ of $V_{k}(S)$ are the centers of the circles through three points of $S$ that enclose exactly $k-1$ or $k-2$ points of $P_k$ and no points of $S\setminus P_k$. If such a circle encloses $k-1$ points, then the vertex is of type I; otherwise it is of type II.
	The interior of an edge $\overline{b_{ij}}$ of a face $f(P_k)$ in $V_k(S)$ consists of exactly those points of $b_{ij}$ that are centers of circles through one point of $P_k$, enclosing the remaining $k-1$ points of $P_k$ and no point of $S\setminus P_k$.
\end{property}

\begin{figure}[h!]	
	\begin{center}
		\includegraphics[width=1\textwidth]{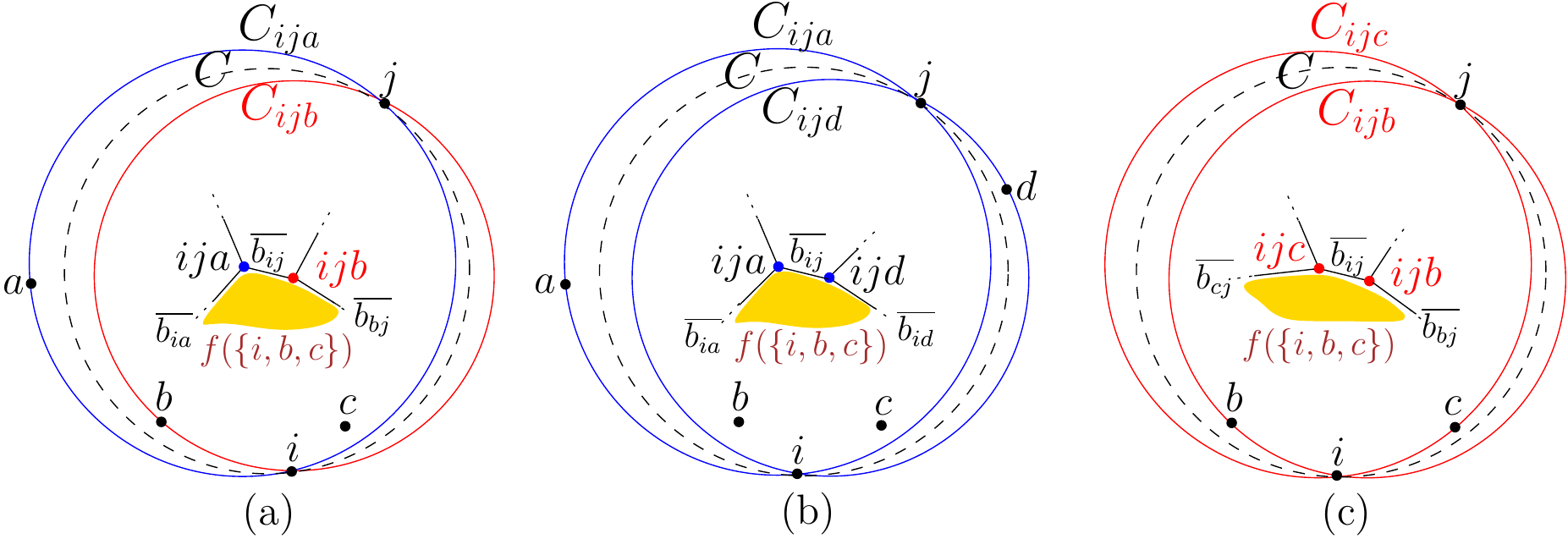}
		\caption{Possibilities for the endpoints of an edge of $V_3(S)$, according to the vertex type. In all three cases, the yellow cell represents a face of $V_3(S)$ defined by $P_3=\{i,b,c\}$. The  circle $C$ (in dashed) goes through $i$ and $j$, and is centered in the interior of $\overline{b_{ij}}$. It always encloses $k-1=2$ points of $P_k$. In case
		(a), $ija$ is of type I (in blue) and $ijb$ is of type II (in red). In case (b), there are two vertices of type I (in blue). And in case (c), there are two vertices of type II (in red). }
		\label{fig:F11-12-22}
	\end{center}
\end{figure}

\begin{proof}
	 Let $ijx$ be a vertex of the face $f(P_k)$. Then $ijx$ is the intersection point of the three bisectors of the points $i,j,x \in S$. Points $i$, $j$ and $x$ are equidistant to the vertex $ijx$.
	In addition, by Property~\ref{property:cell-edge}, at least one of $\{i,j,x\}$, say $i$, belongs to $P_k$, and at least one of them, say $j$, does not belong to $P_k$.
	
	First, assume that $x \notin P_k$. Then
	$ijx$ is of type I because it is the center of a circle through $i\in P_k$ and $j,x\notin P_k$. Furthermore, the $k+2$ nearest points of $S$ to $ijx$ are $P_k\cup \{j,x\}$ by definition of $f(P_k)$. The circle $C_{ijx}$ encloses $P_{k}\setminus\{i\}$, and no points of $S\setminus P_{k}$. See Figure~\ref{fig:F11-12-22} (a) and (b), where $a$ plays the role of $x$.
	
	Second, assume that $x \in P_k$. Then $x \in P_k$, and vertex $ijx$ is of type II. The $k+1$ nearest points of $S$ to $ijx$ are $P_k\cup \{j\}.$
	Then, the circle $C_{ijx}$ encloses $P_{k}\setminus\{i,x\},$ and no points of $S\setminus P_{k}.$ See Figure~\ref{fig:F11-12-22} (a) and (c), where $b$ plays the role of $x$. 
	
	Now consider the edge $\overline{b_{ij}}$ with endpoints $ijx_1$ and $ijx_2$, and let $C$ be a circle through $i$ and $j$ with center in the interior of $\overline{b_{ij}}$. Note that $C$ does not pass through any other point of $S$.
	Note also that $P_k \subset D_{ijx_1}\cap D_{ijx_2}\subset D$ and $D\subset D_{ijx_1}\cup D_{ijx_2}$, where $D$ is the closed disk with boundary $C$. See Figure~\ref{fig:F11-12-22}.
	Hence, $C$ passes through one point $i\in P_k$, and encloses the remaining $k-1$ points of $P_k$ and no point of $S\setminus P_k$.

\end{proof}

\begin{property}\label{property:same-type}
	Let $f^b(P_k)$ be a bounded cell of $V_k(S)$, with $k>1$. Then not all the vertices of $f^b(P_k)$ can be of the same type. 
\end{property}
\begin{proof}

	By Property~\ref{property:cell-edge}, any edge delimiting $f^b(P_k)$ is a segment  $\overline{b_{i_jm_{\ell}}}$ of the bisector between points $i_j \in P_k$ and $m_{\ell}\in S\setminus P_k$.
	
	First, suppose that all vertices of  $f^b(P_k)$ are of type II. Let $i_1,\, i_2\in P_k$ and $m_1,\, m_2\in S\setminus P_k$ such that $\overline{b_{i_1m_1}}$ and $\overline{b_{i_2m_2}}$ are two adjacent edges of $f^b(P_k)$. By Property~\ref{property:ija-ijb}, their common vertex is the center of a circle through three of the points $i_1,i_2,m_1, m_2$. Since the vertex is of type II, we have $m_1=m_2=m$, and the vertex is of the form $i_1i_2m$. 	Hence, the bisectors that delimit $f^b(P_k)$ must be $b_{i_1m},b_{i_2m},\cdots,b_{i_rm}$, see Figure~\ref{fig:FtotsII}.
	By Property~\ref{property:cell-edge}, , for each $j=1,\ldots,r$, point $m$ must lie in the half-plane defined by $b_{i_jm}$  which does not contain $f^b(P_k)$. But this is impossible because $f^b(P_k)$ is a bounded convex region.
	
	\begin{figure}[h!]	
		\begin{center}
			\includegraphics[width=0.5\textwidth]{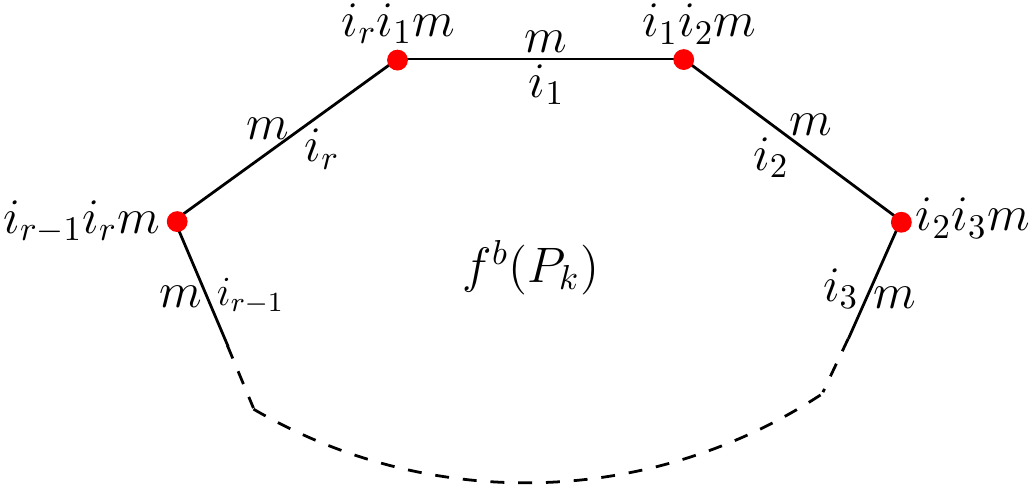}
			\caption{Not all vertices in $f^b(P_k)$ can be of type II.}
			\label{fig:FtotsII}
		\end{center}
	\end{figure}
	
	Suppose now that all vertices are of type I.	Then, by Properties~\ref{property:cell-edge} and~\ref{property:ija-ijb},
	two consecutive vertices of $f^b(P_k)$ are of the form $i\,m_1m_2$ and $i\,m_2m_3$ with $i \in P_k$ and $m_\ell \notin P_k,\, \ell\in\{1,2,3\}$; and the segment connecting them is $\overline{b_{im_2}}$. Consequently, all bisectors delimiting $f^b(P_k)$ are of the form $b_{i\,m_1},b_{i\,m_2},\cdots,b_{i\,m_r}$, with  $i \in P_k$ and $m_\ell \notin P_k,\, \ell\in\{1,\ldots,r\}$, see Figure~\ref{fig:FtotsI}.
	\begin{figure}[h!]	
		\begin{center}
			\includegraphics[width=0.5\textwidth]{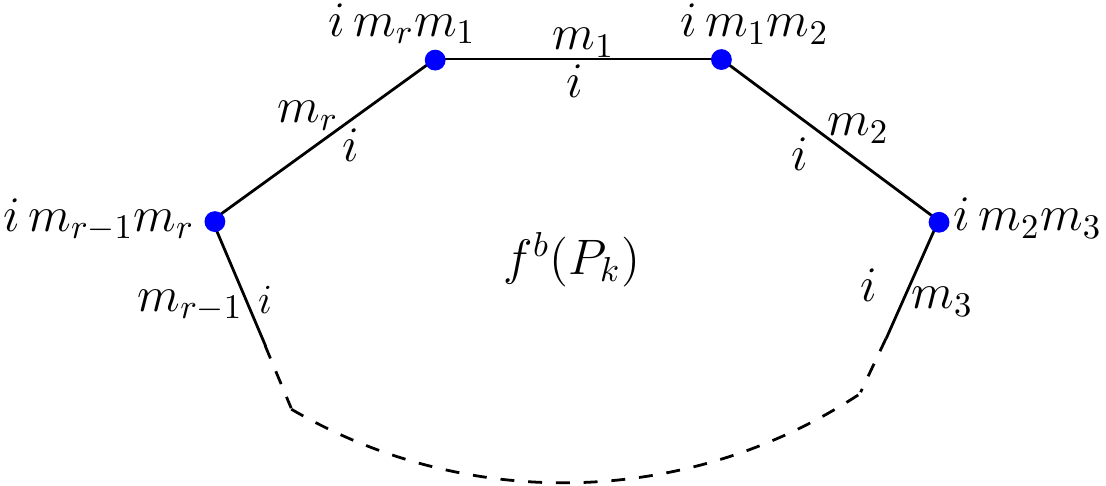}
			\caption{Not all vertices from $V_k(S)$ in $f^b(P_k)$ can be of the type I.}
			\label{fig:FtotsI}
		\end{center}
	\end{figure}
	Then, by Property~\ref{property:cell-edge}, $i$ is inside $f^b(P_k)$.
	For each edge of the boundary of $f^b(P_k)$,  connecting vertices $i\,m_sm_{s+1}$ and $i\,m_{s+1}m_{s+2},\, 1\leq s\leq r$ (with $m_{r+1}=m_{1}$, and $m_{r+2}=m_{2}$), we denote by $\ell_s$ the line through $i$, parallel to this edge. See Figure~\ref{fig:Finter}. On the one hand, the closed half-plane defined by $\ell_s$ wich contains the segment connecting $i\,m_sm_{s+1}$ and $i\,m_{s+1}m_{s+2}$, must also contain $D_{i\,m_sm_{s+1}}\cap D_{i\,m_{s+1}m_{s+2}}$.
	On the other hand, by Property~\ref{property:discs}, all points of $P_k$ must be in the intersection $D_{i\,m_sm_{s+1}}\cap D_{i\,m_{s+1}m_{s+2}},\, \mbox{for}\  1\leq s\leq r$.
	But this would only be possible for $k=1$, since the intersection of the half-planes defined by the lines $\ell_{s}$ through $i$ is the point $i$.
	
	\begin{figure}[h!]	
		\begin{center}
			\includegraphics[width=0.4\textwidth]{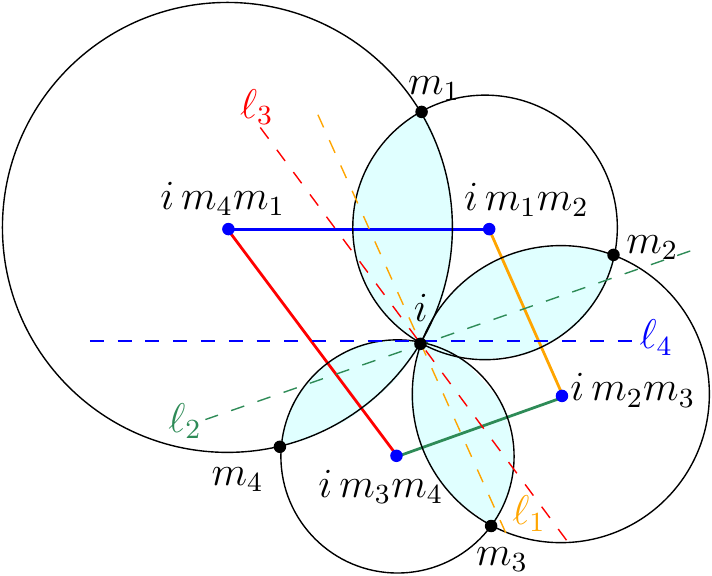}
			\caption{All vertices of $f^b(P_k)$ are of type I. The only point that is in the intersection of all disks is $i$.}
			\label{fig:Finter}
		\end{center}
	\end{figure}
\end{proof}	

Note that unbounded faces of $V_k(S)$ may have only one vertex, and for $k=1$ and $k=n-1$ all vertices of $V_k(S)$ are of the same type.\\

For the next property we also refer to Dehne~\cite{DE84} and Lindenbergh~\cite{L03}. Also, next Figures~\ref{fig:F4-1} and~\ref{fig:F4-2} can be seen as more detailed versions of Figure 1 in~\cite{DE84} and Figure 3 in~\cite{L03}.

\begin{property}\label{property:alternate}
	Let $v$ be a vertex of $V_k(S)$. If $v$ is of type II, then $v$ is also a vertex of $V_{k-1}(S)$. Otherwise, if $v$ is of type I, then $v$ is also a vertex of $V_{k+1}(S)$. Further, the three edges incident to $v$ in $V_k(S)$ alternate with the three edges incident to $v$ in $V_{k-1}(S)$, resp. $V_{k+1}(S)$, in cyclic order around $v$.
\end{property}


\begin{proof}
	Given a circle through three points $a,b,c$ of $S$, if we fix two of them and move slightly the center of the circle along their bisector, the third point becomes either an exterior point or interior point to the circle.
	
	Let $v=abc$ be a vertex of a cell $f(P_k)$ in $V_k(S)$.
	By Property~\ref{property:points-types-edges}, $v$ is either of type I or of type II in $V_k(S)$.
	
	In the case that $v$ is of type I, let  $a\in P_k$ and $b, c\notin P_k$, and consider the circle through $a,b,c$. Suppose we fix $b$ and $c$, and move the center of the circle along $b_{bc}$ towards the interior of $f(P_k)$; see Figure~\ref{fig:F4-1}. By Property \ref{property:cell-edge}, this means moving the center of the circle closer to $a$, so $a$ lies in the interior of the new circle through $b$ and $c$, and the two edges of $f(P_k)$ incident to $v$ are $\overline{b_{ab}}$ and $\overline{b_{ac}}$. Precisely, $a$ lies in the convex region bounded by $b_{ab}$ and $b_{ac}$ that contains the cell $f(P_k)$ of $V_k(S)$. 
	We conclude that the new circle passing through $b$ and $c$ with center on $b_{bc}$ in the interior of $f(P_k)$ contains one point more than the circle through $b$ and $c$ centered in $v$. Since $v$ is of type I in $V_k(S)$ and the circle centered in $v$ encloses $k-1$ points of $S$, the new circle, centered on an interior point of the edge $\overline{b_{bc}}$ incident to $v$, encloses
	$k$ points of $S$. By Property~\ref{property:points-types-edges},
    $\overline{b_{bc}}$ is an edge of $V_{k+1}(S)$. Vertex $v$ is in $V_{k+1}(S)$ because it is an endpoint of $\overline{b_{bc}}$.
	
	We remark that if we move the center of the circle towards the exterior of $f(P_k)$, see Figure~\ref{fig:F4-1} (left), $a$ lies outside of the circle, and $\overline{b_{bc}}$ delimits two cells in $V_k(S)$ defined by $(P_k\setminus \{a\})\cup \{b\}$ and $(P_k\setminus \{a\})\cup \{c\}$, respectively. Here  $\overline{b_{bc}}$ denotes an edge of $V_k(S)$, whereas in the previous paragraph it denoted an edge of $V_{k+1}(S)$ in the same bisector $b_{bc}$. Both edges are incident to $v$ and consecutive along $b_{bc}$.

	\begin{figure}[h!]	
		\includegraphics[width=0.9\textwidth]{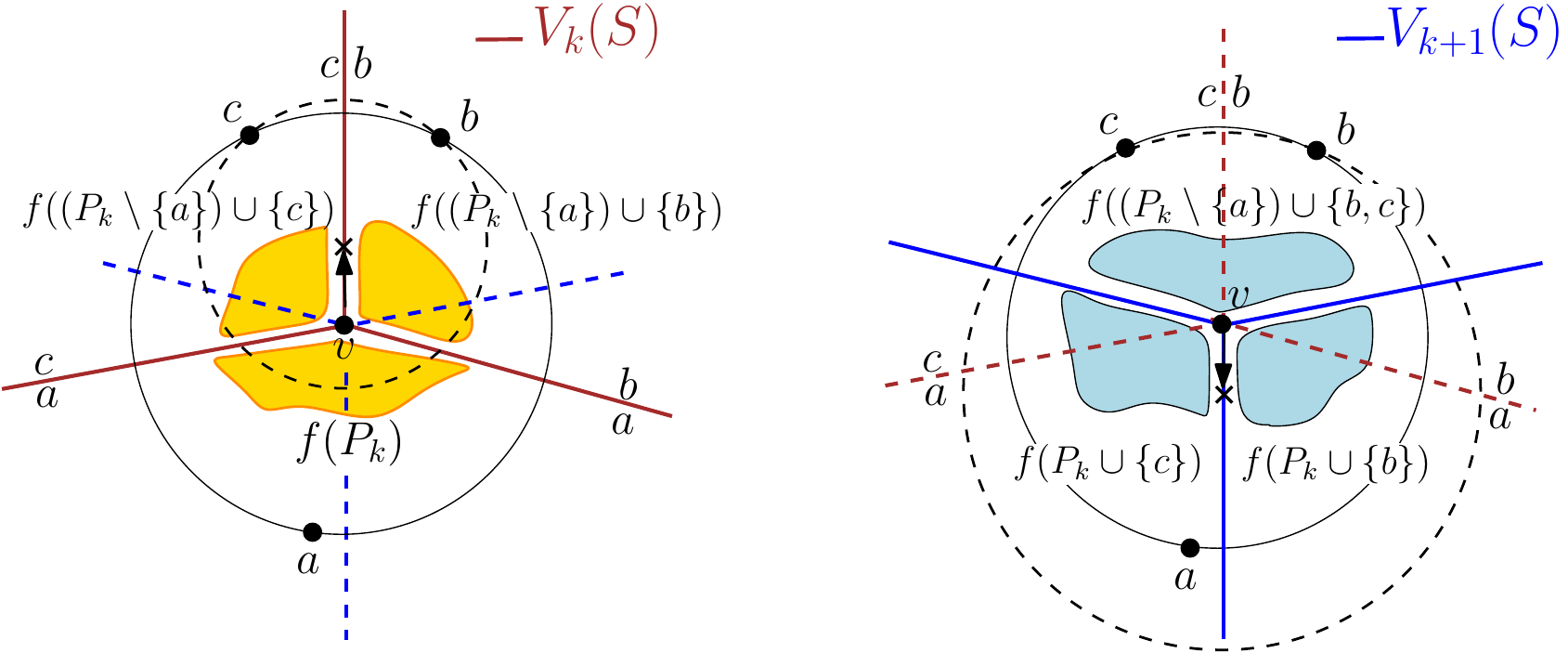}
		\caption{The vertex $v$ is of type I in $V_k(S)$, $a\in P_k$ and $b,c\notin P_k$. Left: The yellow regions represent the three faces of $V_k(S)$ incident to $v$. 
		When $v$ moves outside $f(P_k)$ along the edge $\overline{b_{bc}}$ of $V_k(S)$, point $a$ becomes outside the dashed circle. Right: The blue regions represent the three faces of $V_{k+1}(S)$ incident to $v$. When $v$ moves inside $f(P_k)$ along the edge $\overline{b_{bc}}$ of $V_{k+1}(S)$, point $a$ becomes inside the dashed circle. }
		\label{fig:F4-1}
	\end{figure}
	
	The case that $v$ is of type II is analogous. Let $a,b\in P_k$ and $c\notin P_k$. When we fix $a$ and $b$ and move the center of the circle through $a,b,c$ along $b_{ab}$ towards the interior of the cell $f(P_k)$ of $V_k(S)$, $c$ lies outside the new circle. Then $b_{ab}$ delimits two cells of $V_{k-1}(S)$ defined by $P_k\setminus \{a\}$ and $P_k\setminus \{b\}$, see Figure~\ref{fig:F4-2} (right).
	If the center of the circle through $a,b,c$  moves along $b_{ab}$ towards the exterior of $f(P_k)$, then $c$ lies in the interior of the new circle and $\overline{b_{ab}}$ delimits two cells in $V_{k}(S)$ defined by $(P_k\setminus \{a\})\cup\{c\}$ and $(P_k\setminus \{b\})\cup\{c\}$, respectively. See Figure~\ref{fig:F4-2} (left).
	
	\begin{figure}[h!]	
		\includegraphics[width=0.9\textwidth]{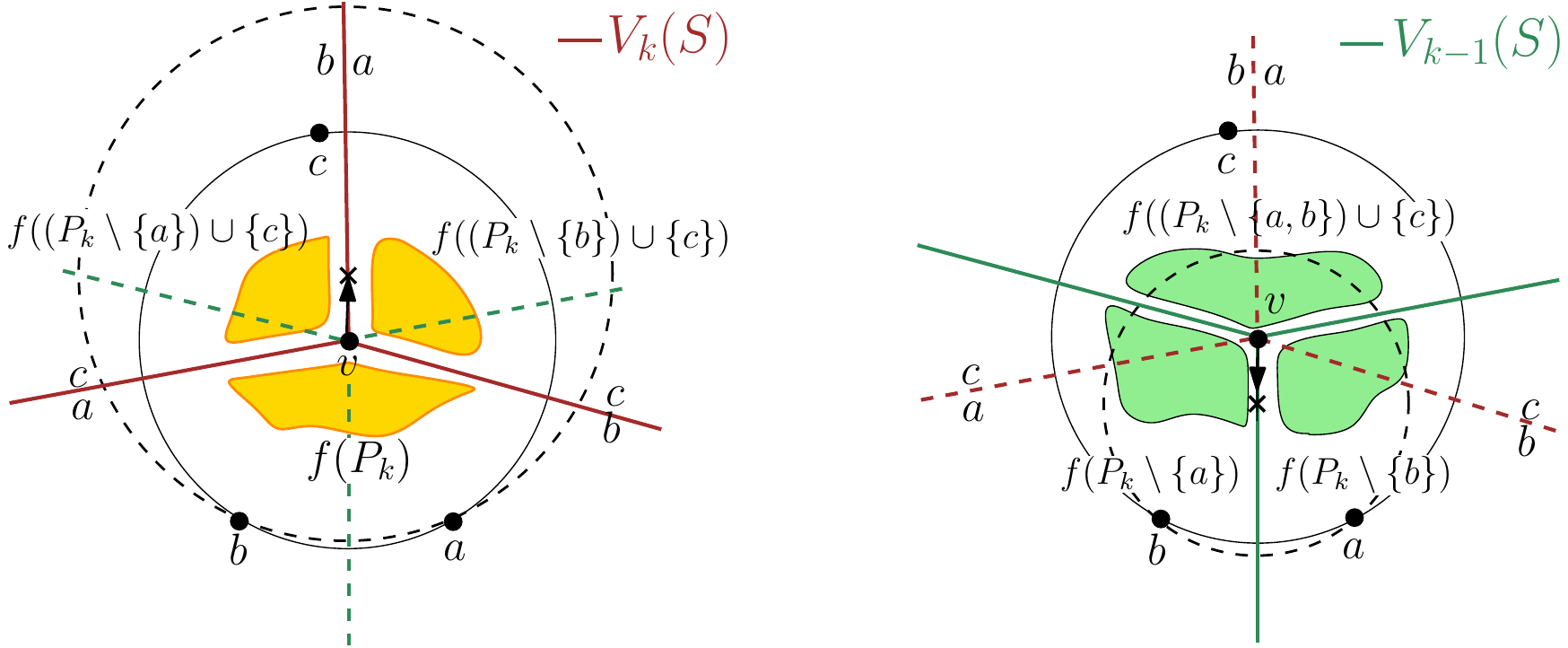}
		\caption{The vertex $v$ is of type II in $V_k(S)$, $a,b\in P_k$ and $c\notin P_k$. Left: The yellow regions represent the three incident faces to $v$ in $V_k(S)$.
	    When $v$ moves outside $f(P_k)$ along the edge $\overline{b_{ab}}$ of $V_k(S)$, point $c$ becomes in the interior of the dashed circle.  
	    Right: The green regions represent the three faces of $V_{k-1}(S)$ incident to~$v$. When $v$ moves inside $f(P_k)$ along the edge $\overline{b_{ab}}$ of $V_{k-1}(S)$, point $c$ becomes outside the dashed circle.  }
		\label{fig:F4-2}
	\end{figure}
	
	Note that every vertex has degree three in $V_k(S)$ and its type is the same for any cell of $V_k(S)$ incident to it. Therefore, applying the previous argument to each cell, the three edges of $V_k(S)$ incident to $v$ alternate with the three edges of $V_{k-1}(S)$ (resp. $V_{k+1}(S)$) incident to $v$. See Figures~\ref{fig:FV234},~\ref{fig:F4-1},~\ref{fig:F4-2}~and~\ref{fig:tree}.
\end{proof}

Each face $f(P_k)$ of $V_k(S)$ is induced by some $P_k\subset S$ with cardinality $\vert P_{k}\vert=k$. Faces of $V_k(S)$ can be divided into subregions which consist of the points of $f(P_k)$ with the same $k+1$ neighbors in $S$. It is known  that this subdivision is achieved by merging appropriately each face $f(P_k)$ with the diagram $V_1(S\setminus P_k)$. It is the basis of the algorithm presented by Lee~\cite{L82} to obtain the Voronoi diagram of order $k+1$.
The following property describes the bisectors that can traverse a cell of $V_k(S)$, corroborating this construction.  

\begin{property}\label{property:edges}
	{Let $f(P_k)$ be a face of $V_k(S)$. Any edge $\overline{b_{ij}}$ of the induced graph of $V_{k+1}(S)$ in $f(P_k)$} has labels $i,j\notin P_k$. And any edge $\overline{b_{ij}}$ of the induced graph of $V_{k-1}(S)$ in $f(P_k)$ has labels $i,j\in P_k$.
\end{property}
\begin{proof}
	Let $\overline{b_{ij}}$ be an edge of the induced graph of $V_{k+1}(S)$ in $f(P_k)$.  By Property~\ref{property:adjacentvoronoi},$\overline{b_{ij}}$ delimits two faces of $V_{k+1}(S)$ whose defining sets of points are $P_k \cup \{i\}$ and $P_k \cup \{j\}$, respectively. We conclude that $i \notin P_k$ and $j \notin P_k$.
	
	Let $\overline{b_{ij}}$ be an edge of the induced graph of $V_{k-1}(S)$ in $f(P_k)$. By Property~\ref{property:adjacentvoronoi}, the edge $\overline{b_{ij}}$ delimits two faces of $V_{k-1}(S)$ whose defining vertices are $P_k\setminus \{i\}$ and $P_k \setminus \{j\}$, respectively. Then $i \in P_k$ and $j \in P_k$.	
\end{proof}
	
\begin{property}\label{property:planar-graph}
	{$V_k(S)\cup V_{k+1}(S)$ is a planar graph for any $1\leq k \leq n-2$.}
\end{property}

\begin{proof}
    First, note that $V_k(S)$, when interpreted as a graph, is a planar graph: Since faces of $V_k(S)$ are convex, if two edges intersect in their interior, then two of their adjacent faces must also intersect in their interior, which contradicts the definition of a Voronoi face.

	Suppose that there is an intersection point $x$ in the interior of two edges,   $\overline{b_{ij}}$ of $V_k(S)$, and $\overline{b_ {\ell m}}$ of $V_{k+1}(S)$. The edge $\overline{b_{ij}}$ is incident to two faces, defined by two sets of $k$ points, $P_k$ and $Q_k$. By Property~\ref{property:adjacentvoronoi}, $Q_k= \left(P_{k} \setminus \{i\}\right) \cup \{j\}.$ Then, the $k+1$ nearest neighbors in $S$ of $x$ are the points in $P_k \cup Q_k$. The edge $\overline{b_{\ell m}}\in V_{k+1}(S)$ traverses both cells $f^b(P_k)$ and $f^b(Q_k)$. Therefore, by Property~\ref{property:edges}, the points $\ell$ and $m$ cannot belong to $P_k\cup Q_k$.
	Moreover, since $x\in\overline{b_ {\ell m}}$, the $k+2$ nearest neighbors of $x$ are the points of $P_k\cup Q_k\cup\{\ell,m\}$.
	This implies that $\left\lbrace \ell, m\right\rbrace \cap \left\lbrace i, j\right\rbrace \neq \emptyset$, which is a contradiction. Then, $V_k(S)\cup V_{k+1}(S)$ is a planar graph. See Figure~\ref{fig:FV234}.
	
	
\end{proof}

\begin{property}\label{property:connected}
	Let $f^b(P_k)$ be a bounded cell of $V_k(S)$. Then, the graph induced by $V_{k-1}(S)$ in $f^b(P_k)$ and the graph induced by $V_{k+1}(S)$ in $f^b(P_k)$ are connected.
\end{property}
\begin{proof}
    Suppose that the graph induced by $V_{k-1}(S)$ in $f^b(P_k)$ is not connected. Let $f^b(P_{k-1})$ be a cell of $V_{k-1}(S)$ that overlaps with $f^b(P_k)$, see Figure~\ref{fig:Fnoconnected} (a).  Note that since the $k-1$ nearest neighbors of any point in $f^b(P_{k-1}) \cap f^b(P_{k})$ are the same, then $P_{k}=P_{k-1}\cup \{m\}$ for some $m\in S$. There are two types of edges in the boundary of $f^b(P_{k-1})\cap f^b(P_{k})$, see Figure~\ref{fig:Fnoconnected} (a):
	\begin{enumerate}
	    \item Edges of $V_{k-1}(S)$, that are in the boundary of $f^b({P_{k-1}})$ and in the interior of $f^b(P_{k})$.
	     By Property~\ref{property:cell-edge}, any edge of $f^b(P_{k-1})\cap f^b(P_k)$ in $V_{k-1}(S)$ is of the form $\overline{b_{xy}},\,$ with $ x\in P_{k-1},\, y\notin P_{k-1}$. By Property~\ref{property:edges}, then $y \in P_{k}$. Since $y \in P_k$ and $y \notin P_{k-1}$, necessarily $y=m$ and $\overline{b_{xy}} = \overline{b_{i_jm}}$ for some $i_j \in P_{k-1}.$ 
	The same argument applies to another edge of this type in another connected component of the graph induced by $V_{k-1}(S)$ in $f^b(P_k)$. Let us denote it by $\overline{b_{i_{\ell}m}}$.
	    \item Edges of $V_k(S)$, that are in the boundary of $f^b(P_{k})$ and in the interior of $f^b({P_{k-1}})$. Let $\overline{b_{xy}}$ be such an edge.
	Suppose $x\in P_k$, then by Property~\ref{property:cell-edge}, $y \notin P_k$. Since $P_{k-1} \subset P_k$, $y\notin P_{k-1}$. In addition, $x \notin P_{k-1}$ by Property~\ref{property:edges}. Then, $x=m$ and $y=a$ for some $a \notin P_k$, so $\overline{b_{xy}} = \overline{b_{ma}}$. The same argument applies to another edge of this type in another connected component of the graph induced by $V_{k-1}(S)$ in $f^b(P_k)$. Let us denote it by $\overline{b_{mb}}$.
	\end{enumerate}
	
	By Property~\ref{property:cell-edge}, $m$ is in the two half-planes defined by $b_{ma}$ and by  $b_{mb}$ that contain $f^b(P_k)$, but does not belong to the two half-planes defined by $b_{i_jm}$ and by $b_{i_{\ell}m}$ that contain $f^b(P_{k-1})$. But this is not possible because the intersection of the respective half-planes containing $m$ is empty.

	\begin{figure}[h!]	
		\begin{center}
			\includegraphics[width=0.8\textwidth]{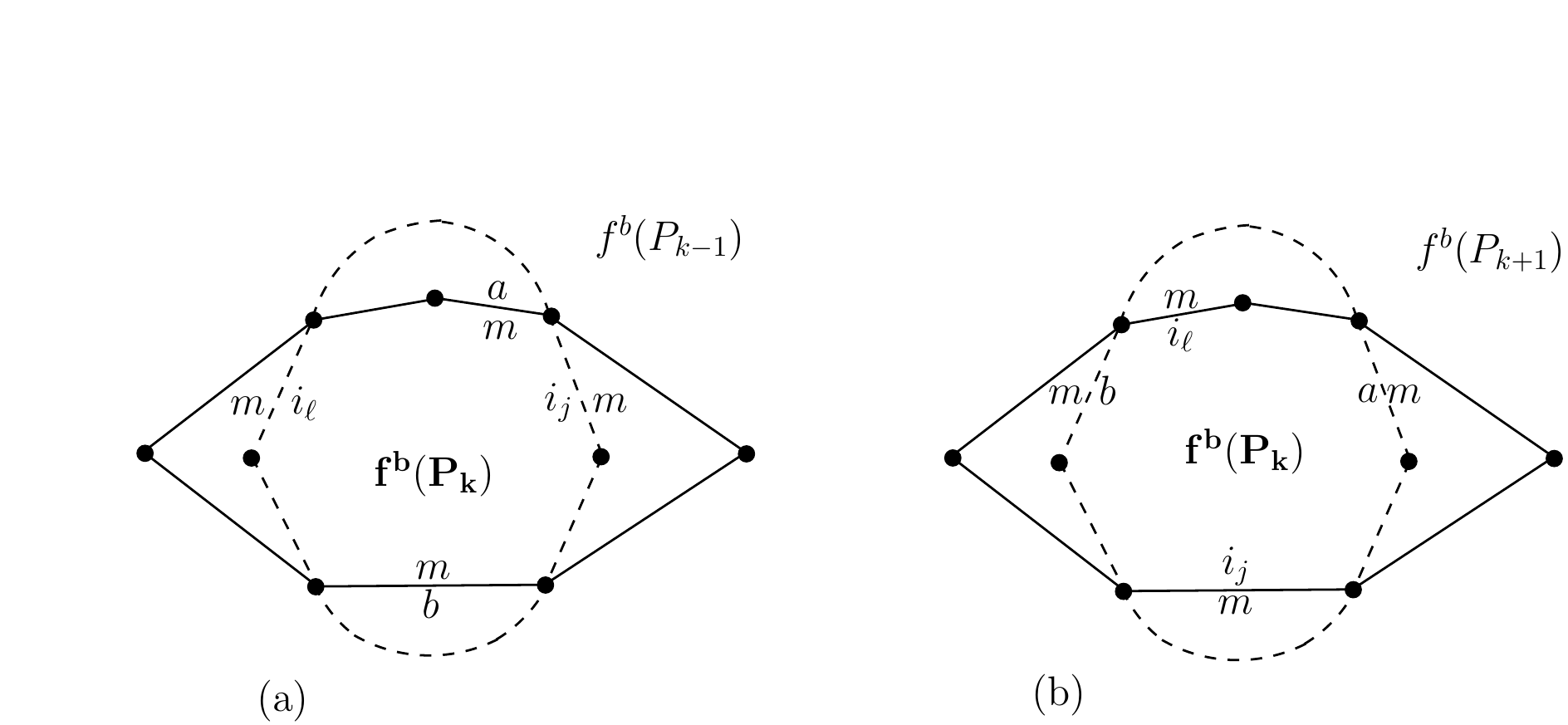}
			\caption{(a) Overlap of the faces $f^b(P_{k})$ from $V_{k}(S)$ and $f^b(P_{k-1})$ from $V_{k-1}(S)$. The face $f^b(P_{k})$ is represented in bold line and 	the face $f^b(P_{k-1})$ is in dashed line. (b) Overlap of the faces $f^b(P_{k})$ from $V_{k}(S)$, in bold line, and $f^b(P_{k+1})$ of $V_{k+1}(S)$, in dashed line.
			The induced subgraphs of $V_{k-1}(S)$ and of $V_{k+1}(S)$ in $f^b(P_k)$ are not connected, which is not possible.}		\label{fig:Fnoconnected}
		\end{center}
	\end{figure}
	
	The same argument serves to prove that the graph induced by $V_{k+1}(S)$ in $f^b(P_k)$ is connected, see Figure~\ref{fig:Fnoconnected} (b).
\end{proof}

Note that Property~\ref{property:connected} can be extended to the graph induced by $V_{k-1}(S)$ in an unbounded face of $V_k(S)$, see Property~\ref{prop:unbdd-forest}, but cannot be extended to the graph induced by $V_{k+1}(S)$.

\begin{property}\label{property:nocycle}
	Let $f(P_k)$ be a cell of $V_k(S)$, with $k>1$. Then the graphs induced by $V_{k-1}(S)$ and $V_{k+1}(S)$ in $f(P_k)$ do not contain cycles.
\end{property}
\begin{proof}
Suppose towards a contradiction that there is a cycle of $V_{k-1}(S)$ contained in $f(P_k)$. Then, since edges of $V_{k-1}(S)$ must alternate with edges of $V_k(S)$ around any common vertex by Property~\ref{property:alternate}, the vertices of this cycle cannot be also vertices of $V_k(S)$. Hence this cycle of $V_{k-1}(S)$ is in the interior of $f^b(P_k)$. For $k=2$ this is a contradiction because all the vertices of $V_1(S)$ are also vertices of $V_2(S)$. Otherwise, for $k >2$ all the vertices of the cycle are vertices of $V_{k-1}(S)\cap V_{k-2}(S)$. In particular, they are all vertices of the same type in a bounded face of $V_{k-1}(S)$, contradicting Property~\ref{property:same-type}. The same argument applies to a face of $V_{k+1}(S)$ inside $f(P_k)$.

	
\end{proof}	

\begin{property}\label{prop:configurations}
  	Let $f^b(P_k)$ and $f^b(P_{k-1})$ be two bounded cells of $V_k(S)$ and $V_{k-1}(S)$, respectively, sharing an interior point. Then $f^b(P_k)$ and $f^b(P_{k-1})$ share exactly two vertices, and neither $f^b(P_k) \subseteq f^b(P_{k-1})$ nor $f^b(P_{k-1}) \subseteq f^b(P_{k})$. That is, among the four configurations of Figure~\ref{fig:Foverlap} only (a) is possible.		
\end{property}
\begin{proof}
	Consider two bounded cells of consecutive higher order Voronoi diagrams, $f^b(P_k)$ and $f^b(P_{k-1})$, having a common interior point.
	The four possible cases are shown in  Figure~\ref{fig:Foverlap}.
	
	Case (b), in which one cell is inside the other and the cells have at least one common vertex $v$, is not possible: By Property~\ref{property:alternate},
	the three edges of $V_k(S)$ incident to $v$ have to alternate with the three edges incident to $v$ from  $V_{k-1}(S)$, in cyclic order around $v$. Case (c) is not possible because the graph induced by $V_{k-1}(S)$ in $f^b(P_k)$ and the graph induced by $V_k(S)$ in $f^b(P_{k-1})$  have to be connected by Property~\ref{property:connected}. Finally, case (d) is not possible because of Property~\ref{property:nocycle}. Only case (a) is possible, that is, the two cells share exactly two vertices, and none of the cells is contained in the other.
	
	\begin{figure}[h!]	
		\begin{center}
			\includegraphics[width=0.8\textwidth]{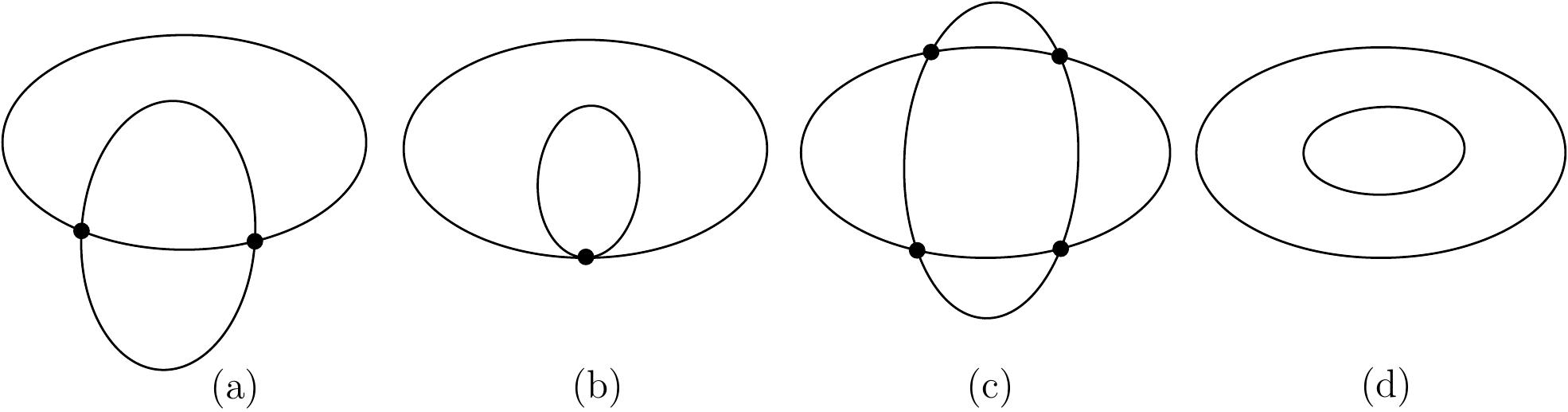}
			\caption{Configurations in which $f^b(P_k)$ and $f^b(P_{k-1})$ have a common interior point. Only case (a) is possible.}	\label{fig:Foverlap}
		\end{center}
	\end{figure}	
\end{proof}

\begin{property}\label{prop:trees}
	Let $f^b(P_k)$ be a bounded cell of $V_k(S)$. For $k>1$, the graph induced by $V_{k-1}(S)$ in $f^b(P_k)$ is a tree with at least two vertices.  Analogously, for $k< n-1$, the graph induced by $V_{k+1}(S)$ in $f^b(P_k)$ is a tree with at least two vertices.	
\end{property}

\begin{proof}
	By Properties~\ref{property:same-type} and~\ref{property:alternate}, the boundary of $f^b(P_k)$ contains vertices of both $V_{k-1}(S)$ and $V_{k+1}(S)$; further the edges of $V_{k+1}(S)$ (analogously, $V_{k-1}(S)$) and of $V_{k}(S)$ incident to a common vertex $v$ alternate in cyclic order. Then, every vertex $v$ of $f^b(P_k)$ is an endpoint of an edge of $V_{k+1}(S)$, or of $V_{k-1}(S)$, lying inside $f^b(P_k)$.
	By the planarity of the graph induced by $V_{k}(S) \cup V_{k+1}(S)$ (or $V_{k}(S) \cup V_{k-1}(S)$) in $f^b(P_k)$, the other endpoint of this edge is also inside $f^b(P_k)$; see Property~\ref{property:planar-graph}. Then, $f^b(P_k)$ contains at least two vertices of $V_{k-1}(S)$ and at least two vertices of $V_{k+1}(S)$. By Properties~\ref{property:connected} and~\ref{property:nocycle}, the graph induced by $V_{k-1}(S)$ (or $V_{k+1}(S)$) in $f^b(P_k)$ has no cycles and is connected. Therefore, it is a tree. See Figure~\ref{fig:tree} for an illustration.
\end{proof}	

\begin{figure}[h!]	
	\begin{center}
		\includegraphics[width=0.3\textwidth]{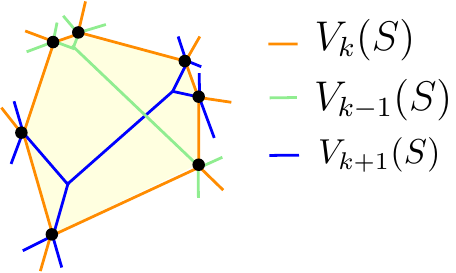}
		\caption{The graph induced by $V_{k-1}(S)$ in a cell $f^b(P_k)$ of $V_k(S)$ is a tree. Analogously, the graph induced by $V_{k+1}(S)$ in $f^b(P_k)$ is a tree.}				
		\label{fig:tree}
	\end{center}
\end{figure}

The following result provides a more precise statement on the number of vertices in each of the two trees of Property~\ref{prop:trees}.

\begin{property}\label{prop:treesize}
    Let $f^b(P_k)$ be a bounded cell of $V_k(S)$ with $\ell$ vertices, for $k>1$. Suppose that $j$ of those vertices also belong to $V_{k-1}(S)$. Then, the graph induced by $V_{k-1}(S)$ in $f^b(P_k)$ is a tree with $2j-2$ vertices, whose $j$ leaves are on the boundary of $f^b(P_k)$. And the graph induced by $V_{k+1}(S)$ inside $f^b(P_k)$ is a tree with $2(\ell-j)-2$ vertices, whose $\ell-j$ leaves are on the boundary of $f^b(P_k)$.
\end{property}
\begin{proof}
	By Property~\ref{prop:trees}, the graph induced by $V_{k-1}(S)$ in $f^b(P_k)$ is a tree $T$ with at least two vertices. Since every vertex of $V_{k-1}(S)$ has degree three, all the leaves of $T$ are on the boundary of $f^b(P_k)$. Moreover, by Property~\ref{property:alternate}, all the vertices of $T$ on the boundary of $f^b(P_k)$ are leaves of $T$. Therefore $T$ has exactly $j$ leaves. If $i$ is the number of interior vertices of $T$, then $T$ has $i+j-1$ edges.	And double counting the number of incidences of vertices and edges of $T$, we get that $j+3i=2(i+j-1)$. This implies that $T$ has $2j-2$ vertices. By Property~\ref{property:alternate}, the $\ell-j$ vertices of $f^b(P_k)$ that are not in $V_{k-1}(S)$ belong to $V_{k+1}(S)$. Hence the tree induced by $V_{k+1}(S)$ in $f^b(P_k)$ has $2(\ell-j)-2$ vertices, whose $\ell-j$ leaves are on the boundary of $f^b(P_k)$.
\end{proof}

    By Properties~\ref{prop:trees} and~\ref{prop:treesize}, the graph induced by $V_{k+1}(S)$ (or $V_{k-1}(S)$) in $f^b(P_k)$ is a tree with its leaves in the boundary of $f^b(P_k)$. The tree gives rise to a subdivision of $f^b(P_k)$ into regions, where each of them is the intersection of $f^b(P_k)$ with a cell of $V_{k+1}(S)$ (resp., $V_{k-1}(S)$).

\begin{figure}[h!]
	\centering
	\subfloat[]{
		\includegraphics[scale=0.4,page=1]{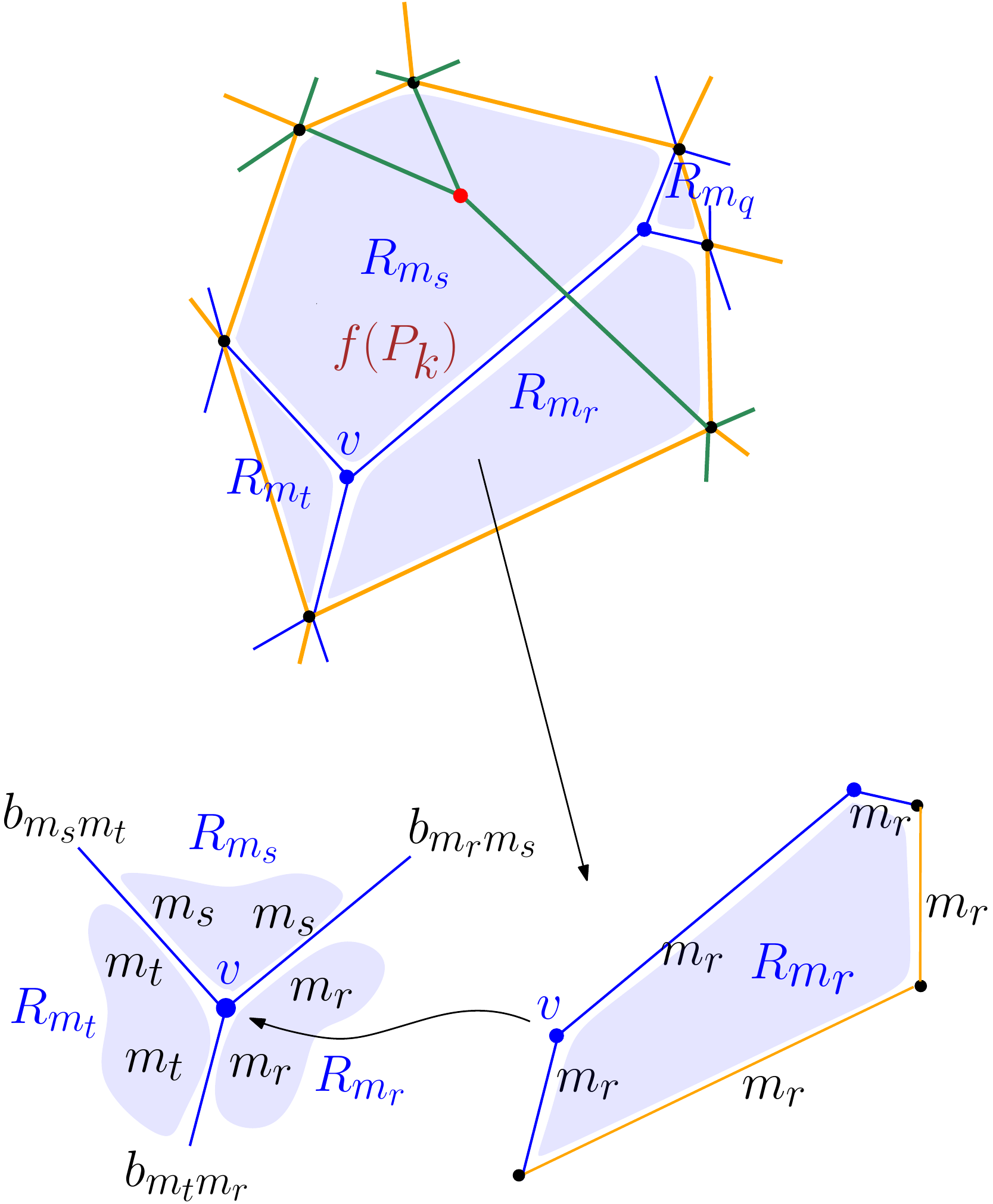}
		\label{fig:k+1}
	}~
	\subfloat[]{
		\includegraphics[scale=0.4,page=2]{Flabels2.pdf}
		\label{fig:k-1}
	}
	\caption{(a) The graph induced by $V_{k+1}(S)$ in $f^b(P_k)$ gives rise to a partition of the cell $f^b(P_k)$  into regions $R_{m_i}$ delimited by bisectors between the $(k+1)$-nearest neighbor $m_i \in S$ of the points of the region, and another point. The cyclic order of the labels of the edges around a vertex $v$ incident with $R_{m_s}$, $R_{m_r}$ and $R_{m_t}$ is $m_r,\,m_s,\,m_s,\,m_t,\,m_t,\,m_r$. (b) The graph induced by $V_{k-1}(S)$ in $f^b(P_k)$ divides the cell $f^b(P_k)$ into regions $R_{i_j}$ delimited by bisectors between the $k$-nearest neighbor $i_j \in S$ of the points of the region, and another point. The cyclic order of the labels of the edges around vertex $v$  incident with $R_{i_r}$, $R_{i_s}$ and $R_{i_t}$, is $i_r,\,i_s,\,i_t,\,i_r,\,i_s,\,i_t$.}
	\label{fig:Labels}
\end{figure}

\begin{property}\label{prop:regions}

    Let $f^b(P_k)$ be a bounded cell of $V_k(S)$, for $k>1$. 
    The edges in the boundary of any region induced by $V_{k+1}(S)$ (resp. $V_{k-1}(S)$) in $f^b(P_k)$
    have  the unique $(k+1)$-nearest neighbor (resp. $k$-nearest neighbor) from $S$ as label inside (outside) the region. Edges in the boundary of $f^b(P_k)$ incident to this region have this label outside (inside) the region. See Figure~\ref{fig:Labels}.
\end{property}
\begin{proof}

    Let $R_{m_i}=f^b(P_k)\cap f^b(P_k\cup\{m_i\})$ be one of the regions induced by $V_{k+1}(S)$ in $f^b(P_k)$. There are two types of edges in the boundary of $R_{m_i}$: edges of $f^b(P_k\cup\{m_i\})$ (edges of the tree) and edges of $f^b(P_k)$, see Figure~\ref{fig:Labels}~(a). On the one hand, the edges of $f^b(P_k\cup\{m_i\})$ that delimit $R_{m_i}$ belong, by Property~\ref{property:cell-edge}, to bisectors between a point from $P_k\cup \{m_i\}$ and another point from $S\setminus (P_k\cup \{m_i\})$. Since these edges are contained in $f^b(P_k)$, by Property~\ref{property:edges}, their labels are not from $P_k$. Then, one of these labels is always $m_i$, the label of the unique $(k+1)$ nearest neighbor from $S$ of the points of $R_{m_i}$. Further, by Property~\ref{property:cell-edge}, the label $m_i$ is inside the region $R_{m_i}$.
    On the other hand, by Property~\ref{property:cell-edge}, the edges from $f^b(P_k)$ that delimit $R_{m_i}$ belong to bisectors between two points, such that one of them is from $P_k$ and the other is from $S\setminus P_k$.
    Since these edges are contained in $f^b(P_k\cup\{m_i\})$, by Property~\ref{property:edges}, their labels are from $P_k\cup\{m_i\}$. Then, one of them is always $m_i$, and this label is outside the region $R_{m_i}$ by Property~\ref{property:cell-edge}.

    Let $R_{i_j}=f^b(P_k)\cap f^b(P_k\setminus\{i_j\})$ be one of the regions induced by $V_{k-1}(S)$ in $f^b(P_k)$. A symmetric argument to the previous case shows that all the edges of $R_{i_j}$ have label $i_j$, where $i_j$ is the unique $k$-nearest neighbor from $S$ of the points of $R_{i_j}$. Furthermore, this label is inside $R_{ij}$ for the edges in the boundary of $f^b(P_k)$ and outside for the others.
  
\end{proof}	

\begin{property}\label{prop:labels}
	Let $v$ be a vertex of $V_{k+1}(S)$ or $V_{k-1}(S)$ in the interior of $f^b(P_k)$.
	If $v$ is from $V_{k+1}(S)$, let $m_r,m_s,m_t$ be the points of $S\setminus P_k$ that define the incident bisectors with $v$.	
	Then the cyclic order of the labels of the edges around $v$ is $m_r,\,m_s,\,m_s,\,m_t,\,m_t,\,m_r$.
	If $v$ is from $V_{k-1}(S)$, let $i_r,i_s,i_t$ be the points of $P_k$ defining the incident bisectors with $v$. Then, the cyclic order of the labels of the edges around $v$ is  $i_r,\,i_s,\, i_t,\,i_r,\,i_s,\,i_t$. 
\end{property}
\begin{proof}
    Let $v$ be a vertex from $V_{k+1}(S)$ in the interior of $f^b(P_k)$. By  Property~\ref{property:edges}, the incident bisectors with $v$ are bisectors between pairs of points of $S\setminus P_k$. These points are $m_r,m_s,m_t$. By Property~\ref{prop:regions},
    all the edges delimiting a region induced by $V_{k+1}(S)$ in $f^b(P_k)$ belong to the bisector between two points, where one of them is the $(k+1)$-closest neighbor from $S$ of the points of the region. Let $R_{m_i}=f^b(P_k)\cap f^b(P_k\cup\{m_i\})$ be the region where the $(k+1)$-closest neighbor from $S$ is $m_i$. Therefore, $v$ is the vertex incident to the regions $R_{m_r},R_{m_s}$ and $R_{m_t}$. Since the two edges incident to $R_{m_i}$ have label $m_i$ in the interior of $R_{m_i}$ (for $i\in\{r,s,t\}$), the cyclic order of the labels of the edges around $v$ is: $m_r,\,m_s,\,m_s,\,m_t,\,m_t,\,m_r$. See Figure~\ref{fig:Labels} (a).

    Now suppose that $v$ is from $V_{k-1}(S)$. By  Property~\ref{property:edges}, the incident bisectors with $v$ are bisectors between pairs of points of $P_k$. These points are $i_r,i_s,i_t$. 
    By Property~\ref{prop:regions},
    every edge delimiting a region induced by $V_{k-1}(S)$ in $f^b(P_k)$ belongs to the bisector between two points; where one of them is the $k$-closest neighbor from $S$ of the points of the region. Let $R_{i_j}=f^b(P_k)\cap f^b(P_k\setminus\{i_j\})$ be the region where the $k$-closest neighbor from $S$ is $i_j$. Therefore, $v$ is the vertex common to the regions $R_{i_r},R_{i_s}$ and $R_{i_t}$. 
    Since the two edges incident to $R_{i_j}$ have label $i_j$ in the exterior of $R_{i_j}$ (for $j\in\{r,s,t\}$), the cyclic order of the labels of the edges around $v$ is: $i_r,\,i_s,\, i_t,\,i_r,\,i_s,\,i_t$. See Figure~\ref{fig:Labels} (b). 
\end{proof}

We restate Property~\ref{prop:labels} in terms of vertices of type I and II. 

\begin{property}[\bf{Vertex rule}]\label{prop:cyclic-order}
    Let $v$ be a vertex of $V_k(S)$ and let $i, j,\ell\in S$ be the points that define the three incident bisectors with $v$. If $v$ is of type I, the cyclic order of the labels around $v$ is $i, i, j, j, \ell, \ell$. And if $v$ is of type II, the cyclic order of the labels around $v$ is $i, j, \ell, i, j, \ell$.
\end{property}
\begin{proof}
    Assume first that $v$ is a vertex of type I of $V_k(S)$ for $k>1$. Then, by Property~\ref{property:alternate}, $v$ is not a vertex of $V_{k-1}(S)$, and therefore $v$ is in the interior of a face of $V_{k-1}(S)$. By Property~\ref{prop:labels}, the cyclic order of the labels around $v$ is $i, i, j, j, \ell, \ell$. Now assume that $v$ is a vertex of type II of $V_k(S)$. Then $v$ is in the interior of a face of $V_{k+1}(S)$. By Property~\ref{prop:labels}, the cyclic order of the labels around $v$ is $i, j, \ell, i, j, \ell$.
    Note that a vertex of type I in $V_k(S)$ is a vertex of type II in $V_{k+1}(S)$. Hence, Property~\ref{prop:cyclic-order} is also valid for $k=1$.
\end{proof}

\begin{figure}[ht]
	\centering
	\includegraphics[scale=0.7]{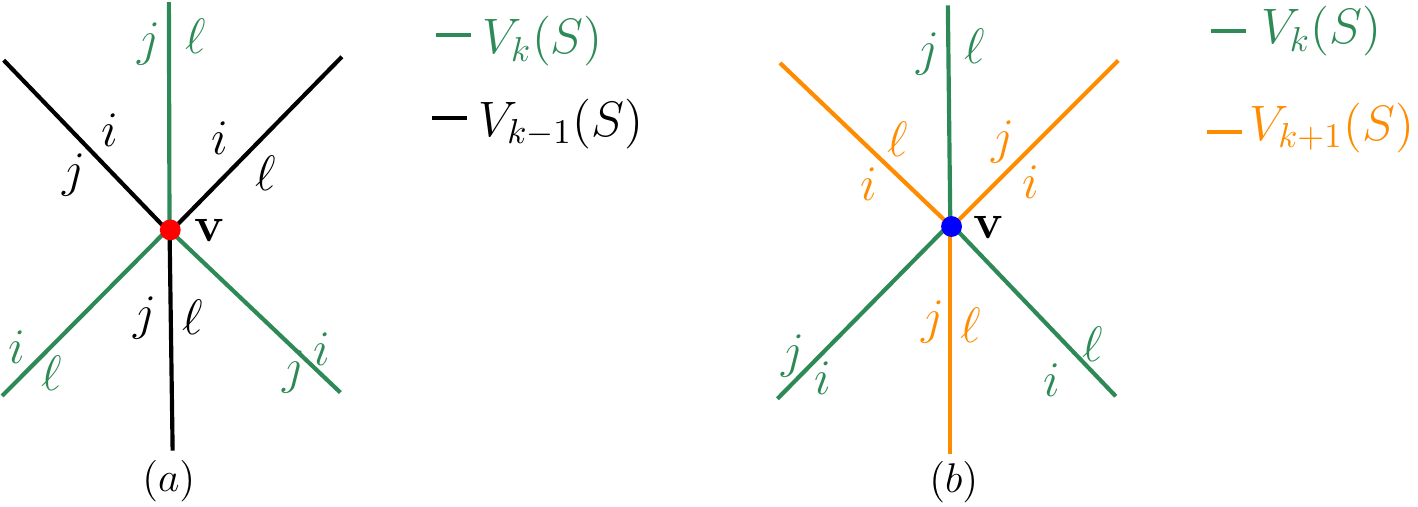}
	\caption{(a) Vertex $v$ is of type II in $V_k(S)$. The edge labeling around $v$ in $V_k(S) \cup V_{k-1}(S)$. (b) Vertex $v$ is of type I in $V_k(S)$. The edge labeling around $v$ in $V_k(S) \cup V_{k+1}(S)$.}
	\label{fig:vertex-label}
\end{figure}

Figure~\ref{fig:vertex-label} illustrates Property~\ref{prop:cyclic-order} for two consecutive Voronoi diagrams. 

\begin{property}\label{prop:trianglefree}
	For $k>1$, every bounded face $f^b(P_k)$ of $V_k(S)$ contains at least two and at most $n-k$ vertices of type I and at least two and at most $k$ vertices of type II. In particular, $V_k(S)$ does not contain triangles.
\end{property}
\begin{proof}
The lower bounds on the numbers of vertices are a consequence of Properties~\ref{prop:trees} and~\ref{prop:treesize}. In particular, they imply that for $k > 1$, $V_k(S)$
does not contain triangles. 

	
	For the upper bound of vertices of type II consider the subdivision defined by $V_{k-1}(S)$ inside $f^b(P_k)$ given by Property~\ref{prop:regions}. Since each region of the subdivision is defined by the points of the plane with the same set of $k$ nearest neighbors but with a different $k$-nearest neighbor, there are at most $k$ such regions (one for each $k$-nearest neighbor). Each such region is incident to exactly two vertices in the boundary of $f^b(P_k)$, and each such vertex is incident to two regions through an edge of $V_{k-1}(S)$. Therefore $f^b(P_k)$ has at most $k$ vertices that also belong to $V_{k-1}(S)$. By Property~\ref{property:alternate}, these are the vertices of type II of $f^b(P_k)$. 
	For the vertices of type I, consider the subdivision into regions defined by $V_{k+1}(S)$ inside $f^b(P_k)$. Each of these regions is defined by the same set of $k$ nearest neighbors $P_k$ and a different $(k+1)$-nearest neighbor in $S\setminus P_k$. Therefore, the subdivision contains at most $n-k$ regions and $f^b(P_k)$ contains at most $n-k$ vertices of type I. 
\end{proof}

\begin{figure}[h!]
		\centering
		\includegraphics[scale=0.7]{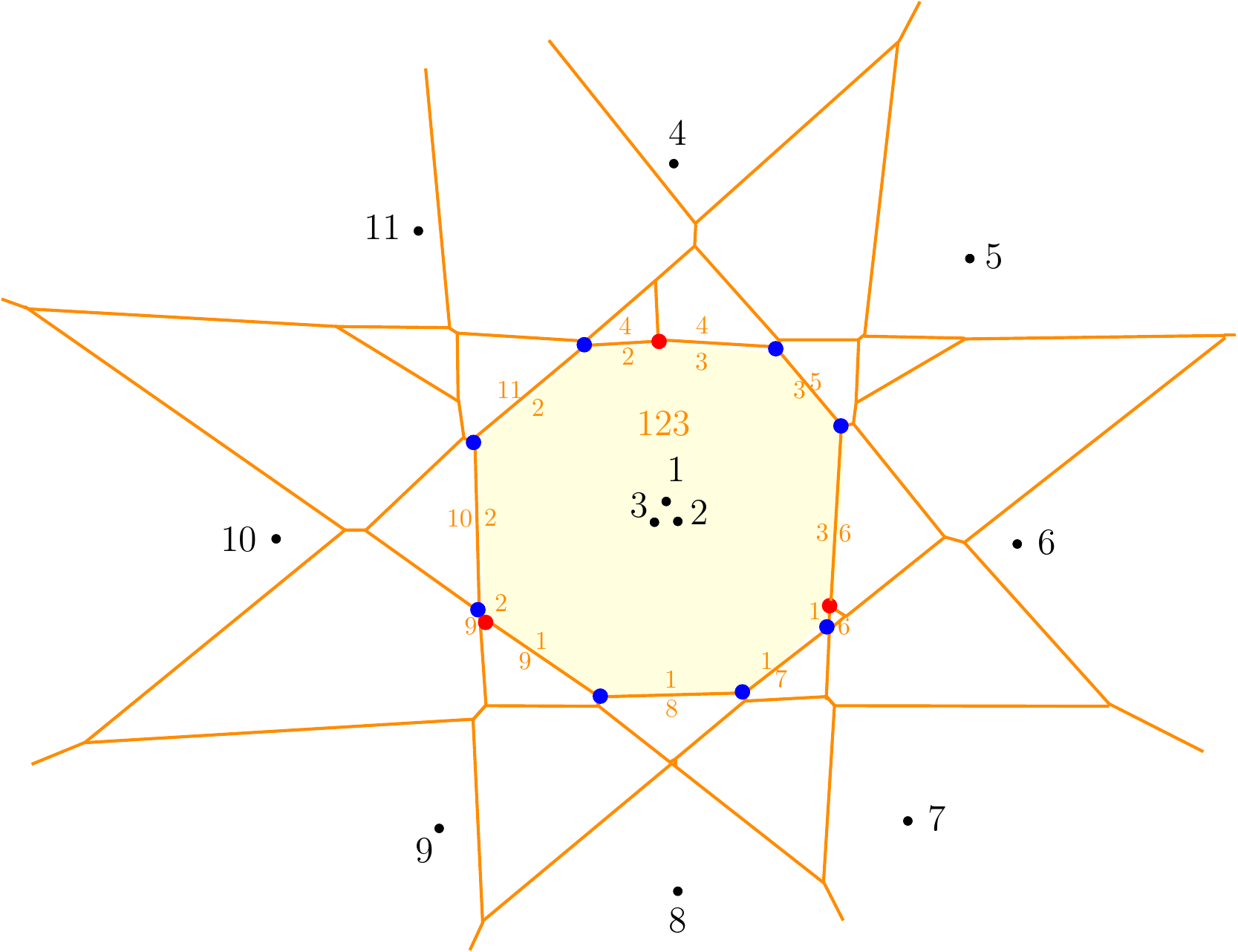}
		\caption{A cell of $V_3(S)$ for a set of $n=11$ points, where both the number of vertices of type I ($n-k=8$) and type II ($k=3$) are maximum.}\label{fig:max-I-II}	
	\end{figure}

The lower bounds in Property~\ref{prop:trianglefree} were already given
in~\cite{L82}. It was also proved in~\cite{MRT19} with a
different method that $V_2(S)$ does not contain triangles. 
As for the upper bounds in Property~\ref{prop:trianglefree},  there exist point sets $S$ such that some face of $V_k(S)$ has exactly $n-k$ vertices of type I and $k$ vertices of type II, see Figure~\ref{fig:max-I-II}.

Note that Property~\ref{prop:trianglefree} implies that every bounded face of $V_2(S)$ contains exactly two vertices of type~II, and every bounded face of $V_{n-2}(S)$ contains exactly two vertices of type~I.

\begin{figure}[h!]
    \centering
    \subfloat[]{
    	\includegraphics[scale=0.4,page=1]{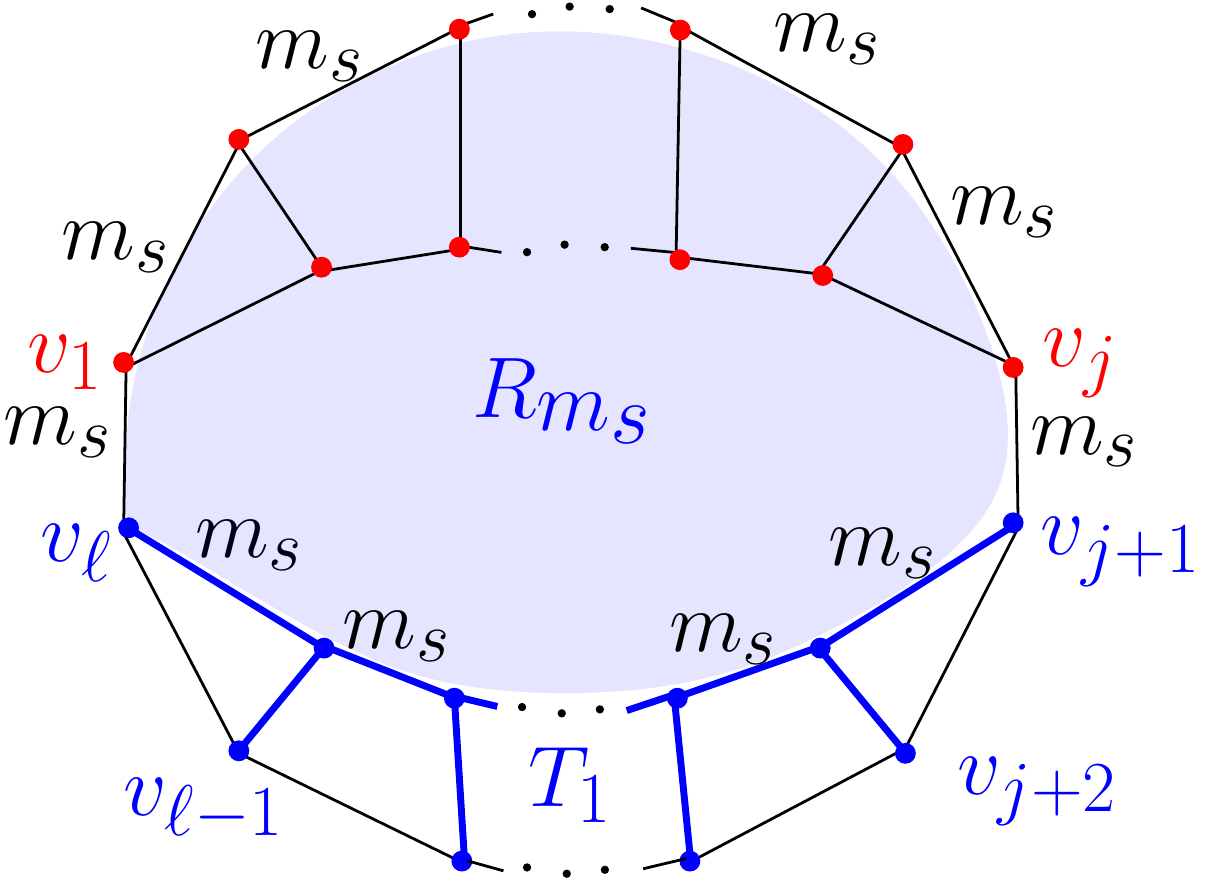}
    }~
    \subfloat[]{
    	\includegraphics[scale=0.4,page=2]{F2-chains.pdf}
    }
    \caption{(a) The edges $v_jv_{j+1}$ and $v_{\ell}v_1$ belong to the bisector between $m_s$ and another point of $P_k$. (b) The edges $v_jv_{j+1}$ and $v_{\ell}v_1$ belong to the bisector between $i_r$ and another point of $S\setminus P_k$.
    }
    \label{fig:2-chains}
\end{figure}

\begin{property}\label{prop:no-two-chains}
	Let $f^b(P_k)$ be a bounded cell of $V_k(S)$, for $k>1$. 
	Then, on its boundary not all vertices of the same type are consecutive. Equivalently, inside $f^b(P_k)$, there is an edge of $V_{k+1}(S)$ that crosses an edge of $V_{k-1}(S)$.
\end{property}
\begin{proof}
	By Property ~\ref{prop:treesize}, the boundary of $f^b(P_k)$ has vertices of type I and II, which are the leaves of the trees induced in $f^b(P_k)$ by $V_{k+1}(S)$ and $V_{k-1}(S)$, respectively. Let $T_1$ and $T_2$ be the corresponding trees. 
	
	Suppose that all the vertices of the same type are consecutive in the ordered list of the boundary vertices of $f^b(P_k)$.
	Let $v_1,\cdots, v_j$ be the vertices of type II and $v_{j+1},\cdots, v_{\ell}$ the vertices of type I of $f^b(P_k)$.
	By Property~\ref{prop:regions}, $T_1$ gives rise to a subdivision of the cell $f^b(P_k)$ into regions which differ in their $(k+1)$-nearest neighbor from $S$. Analogously, $T_2$,
	gives rise to a subdivision into regions which differ in their $k$-nearest neighbor from $S$.
	
	Then, there is a path between $v_{j+1}$ and $v_{\ell}$ in $T_1$, that delimits a region whose $(k+1)$-nearest neighbor from $S$ is a point of $S\setminus P_k$. Let $m_s$ be this point and $R_{m_s}$ the corresponding region, see Figure~\ref{fig:2-chains} (a). Then, by Property~\ref{prop:regions}, the edges $\overline{v_1v_{\ell}}$ and $\overline{v_jv_{j+1}}$ have label $m_s$ outside $R_{m_s}$.
	Analogously, there is a path between $v_1$ and $v_j$ in $T_2$, that delimits a region whose $k$-nearest neighbor from $S$ is a point of $P_k$. Let $i_r$ be this point and $R_{i_r}$ the corresponding region, see Figure~\ref{fig:2-chains} (b).  Then, by Property~\ref{prop:regions}, the edges $\overline{v_1v_{\ell}}$ and $\overline{v_jv_{j+1}}$ have label $i_r$ inside $R_{i_r}$.
	
	Therefore, the edges $v_1v_{\ell}$ and $v_jv_{j+1}$ should belong to the same bisector, $b_{m_si_r}$, which is impossible.
\end{proof}

\begin{property}[\bf{Face rule}]\label{prop:consecutive_labels}
In each face of $V_k(S)$, the edges that have the same label $i$ are consecutive, and these labels $i$ are either all in the interior of the face, or are all in the exterior of the face.
\end{property}

\begin{proof}
By Property~\ref{prop:regions}, a face $f^b(P_k)$ of $V_k(S)$ can be constructed by joining regions of at most $k$ faces of $V_{k-1}(S)$ as described in Figure~\ref{fig:Labels}. Each of these regions is the set of points in the plane with the same set of $k$ nearest neighbors $P_k$, but with a different $k$-nearest neighbor. That is, each point in $P_k$ determines at most one region of $f^b(P_k)$. According to Property~\ref{prop:regions}, if the $k$-nearest neighbor of a region is the point $i$, then the boundary edges of the region that belong to $V_k(S)$ have the interior label $i$. Note that any other region in the subdivision of $f^b(P_k)$ cannot have label $i$ as an interior boundary label.
Since each boundary edge of $f^b(P_k)$ must be a boundary edge of one such region,
then all the boundary edges of $f^b(P_k)$ with interior label $i$ are precisely the edges of the region of $f^b(P_k)$ such that $i$ is the $k$-nearest neighbor for this region. Therefore they are consecutive. The graph of $V_{k-1}(S)$ inside $f^b(P_k)$ is a tree, then there is a single interval of consecutive edges with the same label $i$.

A similar argument is used to show that the edges with label $i$ in the exterior of $f^b(P_k)$ are consecutive. Consider the subdivision in regions of $V_{k+1}(S)$ inside $f^b(P_k)$ given in Property~\ref{prop:regions}. Each of these regions is the set of points with the same set of $k$ nearest neighbors $P_k$, but with a different $(k+1)$-nearest neighbor. If the $(k+1)$-neighbor of a region is the point $i$, then all edges of the region in the boundary of $f^b(P_k)$ (edges of $V_k(S)$) have the exterior label $i$; see Figure~\ref{fig:Labels}. 
Therefore, different regions have a different exterior label of the edges of $f^b(P_k)$. This implies that all the edges of $f^b(P_k)$ with the same exterior label $i$ are consecutive.
Note that by Property~\ref{property:adjacentvoronoi}, the consecutive labels $i$ have to be either all in the interior or all in the exterior of $f^b(P_k)$.
\end{proof}

\section{Properties about unbounded faces of $V_k(S)$}\label{sec:unbounded}

In the previous sections, we have studied a number of properties concerning the bounded faces of $V_k(S)$. In the current section, we consider if these properties can be adapted to the case of unbounded faces, namely Properties~\ref{property:same-type},  \ref{prop:configurations}, \ref{prop:trees}, \ref{prop:regions} and \ref{prop:no-two-chains}.

In the first place, the analogous of Property~\ref{property:same-type} is not true for unbounded faces of $V_k(S)$ because they may have only one vertex.

\begin{figure}[h!]
	\centering
	\includegraphics[scale=0.6]{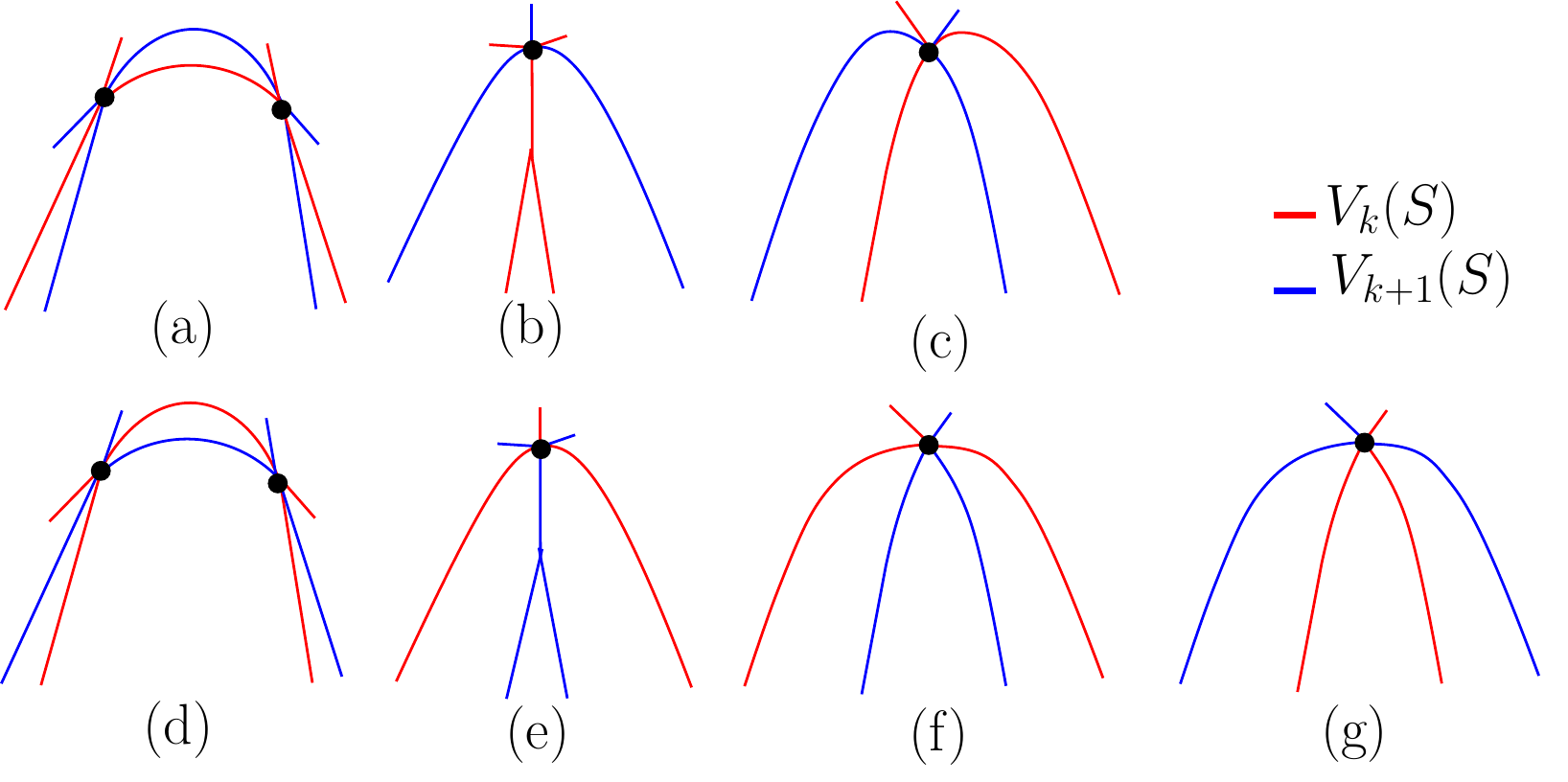}
	\caption{Simplified configurations with unbounded cells: (a), (b) and (c) are possible; (d), (e), (f) and (g) are not possible.}
	\label{fig:unbdd-cells}
\end{figure}

Similar to Property~\ref{prop:configurations}, we have the following characterization of unbounded faces:

\begin{property}\label{prop:conf-unbonded}
Let $f^{\infty}(P_k)$ and $f^{\infty}(P_{k+1})$ be two unbounded cells of $V_k(S)$ and  $V_{k+1}(S)$, respectively, that intersect. Then, among the six configurations represented in Figure~\ref{fig:unbdd-cells}, only three of them are possible: (a), (b) and~(c).		
\end{property}

\begin{proof}
First we note that any unbounded face of $V_k(S)$ always intersects some unbounded face of $V_{k+1}(S)$.
Consider two such unbounded faces  $f^\infty(P_k)$ and $f^\infty(Q_{k+1})$ of consecutive Voronoi diagrams $V_k(S)$ and $V_{k+1}(S)$ with some common vertex.
Figure~\ref{fig:unbdd-cells} shows six configurations which correspond to the following three cases: 1) The two faces $f^\infty(P_k)$ and $f^\infty(Q_{k+1})$ share at least two vertices and (a) the graph induced by $V_{k+1}(S)$ inside $f^\infty(P_k)$ is not connected or (d) the graph induced by $V_{k}(S)$ inside $f^\infty(Q_{k+1})$ is not connected.
 2) The two faces $f^\infty(P_k)$ and $f^\infty(Q_{k+1})$ do not have common vertices and either (b) $f^\infty(P_k) \subset f^\infty(Q_{k+1})$ or (e) $f^\infty(Q_{k+1}) \subset f^\infty(P_{k})$. 3) Both cells share exactly one vertex: (c), (f) and (g). 

The options (d), (e), (f) and (g) are not possible:
In option (d) the graph induced by $V_k(S)$ in $f^{\infty}(Q_{k+1})$ is not a tree contradicting Property~\ref{prop:unbdd-forest}.
In option (e) $f^{\infty}(Q_{k+1})$ is in the interior of $f^{\infty}(P_k)$. Then all the vertices of $f^{\infty}(Q_{k+1})$ are of the same type. This implies that all the edges in the boundary of $f^{\infty}(Q_{k+1})$ have the same interior label, say $i$. Therefore, point $i$ is inside $f^{\infty}(Q_{k+1}) \subset f^{\infty}(P_k)$ and $i\in P_k$. This is a contradiction with  Property~\ref{property:edges}. Finally, the options (f) and (g) are not possible because the vertex rule is not fulfilled.
Examples show that the remaining three configurations (a), (b) and (c) are possible.
\end{proof}

In Property~\ref{prop:trees} we have proved that the graphs induced by $V_{k+1}(S)$ and $V_{k-1}(S)$ inside a bounded face of $V_k(S)$ are trees. In Property~\ref{prop:unbdd-forest} we show that the analogous statement for the graph induced by $V_{k-1}(S)$ inside an unbounded face of $V_k(S)$ is also true. Then, we argue that this is not the case for the graph induced by $V_{k+1}(S)$.

\begin{property}\label{prop:unbdd-forest}
    Let $f^{\infty}(P_k)$ be an unbounded face of $V_k(S)$. Then the subgraph induced by $V_{k-1}(S)$ inside $f^{\infty}(P_k)$ is a tree.
\end{property}

\begin{proof}
First, observe that $f^\infty(P_k)$ cannot contain a cycle of $V_{k-1}(S)$ because such a cycle would be a bounded face such that all its vertices are of the same type in $V_{k-1}(S)$, contradicting Property~\ref{property:same-type}.

    Assume, towards a contradiction, that the subgraph induced by $V_{k-1}(S)$ inside $f^{\infty}(P_k)$ is a forest with more than one connected component. Then all the leaves of this forest must be type II vertices in $V_k(S)$, that is, vertices of $V_k(S)$ that are also vertices of $V_{k-1}(S)$. Consider two such leaves $v_1$ and $v_2$ that belong to different trees $T_1$ and $T_2$. Choose $v_1$ and $v_2$ such that all the vertices in the boundary of $f^{\infty}(P_k)$ between $v_1$ and $v_2$ are vertices of type I, or such that $\overline{v_1v_2}$ is an edge of $V_k(S)$. Then, by Property~\ref{prop:cyclic-order}, all the boundary edges of $f^{\infty}(P_k)$ between $v_1$ and $v_2$ have the same label, say $j$, in the interior of $f^{\infty}(P_k)$, see Figure~\ref{fig:unbdd-forest}. The defining points of vertices $v_1$ and $v_2$ are denoted with $i,j,l,m, x$ according to Figure~\ref{fig:unbdd-forest}.
    Since $v_1$ and $v_2$ are vertices of type II in $V_k(S)$, the cyclic order of the labels of the edges of $V_k(S)$ around $v_1$ and $v_2$ is $ij\ell ij\ell$ and $jxmjxm$, respectively.
    On the other hand, $v_1$ and $v_2$ are vertices of type I in $V_{k-1}(S)$. This implies that the cyclic order of the labels of the edges of $V_{k-1}(S)$ around $v_1$ and $v_2$ is $iijj\ell\ell$ and $jjmmxx$, respectively. 
    The position of label $j$ in the edges of $V_{k}(S)$ also induces a position for the label $j$ in the edges of $V_{k-1}(S)$ around $v_1$ and $v_2$ by Property~\ref{property:alternate}, see Figure~\ref{fig:unbdd-forest}. Then point $j$ must be in the interior of the two yellow cones in Figure~\ref{fig:unbdd-forest}, with apex $v_1$ and $v_2$, respectively. Since $T_1$ and $T_2$ are two different connected components, the two cones cannot intersect. Therefore we get a contradiction for the position of point $j$.

 \begin{figure}[h!]
    	\centering
    	\includegraphics[scale=0.75]{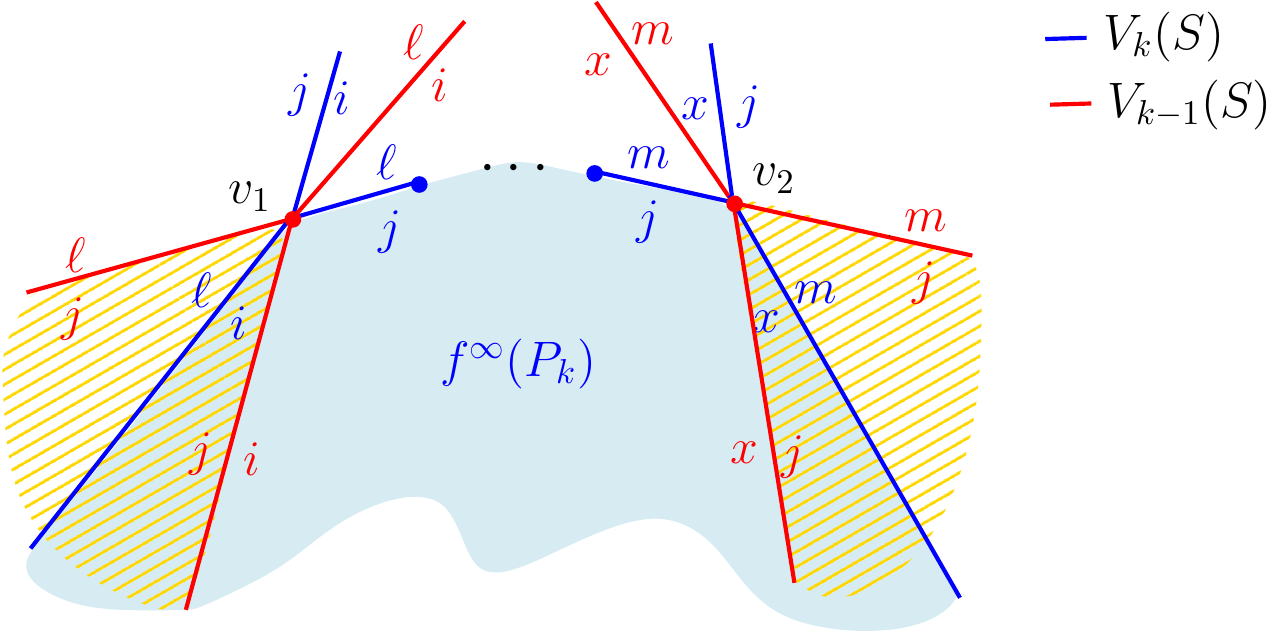}
    	\caption{Illustration of the proof of Property~\ref{prop:unbdd-forest}.}
    	\label{fig:unbdd-forest}
    \end{figure}
\end{proof}

\begin{figure}[h!]
	\centering
	\includegraphics[scale=0.6, page =1]{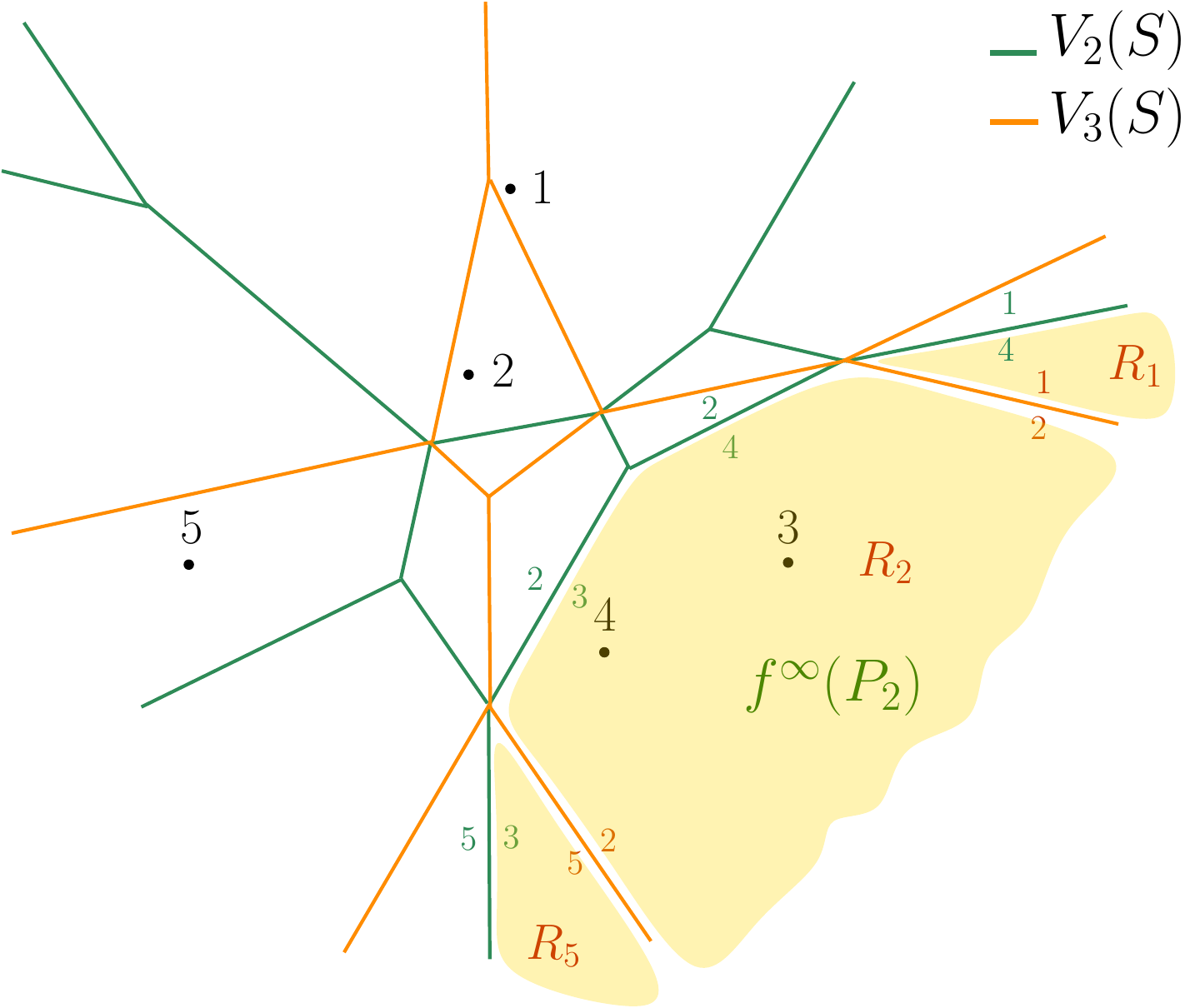}
	\caption{The graph induced by $V_{3}(S)$ inside an unbounded face $f^\infty(P_2)$ of $V_2(S)$ can be a forest with more than one connected component.}
	\label{fig:unbdd-face}
\end{figure}

The analogous statement of Property~\ref{prop:trees} for unbounded faces of $V_{k+1}(S)$ is not true. That is, the graph induced by $V_{k+1}(S)$ inside an unbounded face of $V_k(S)$ can be the empty graph with no vertices or a forest with more than one connected component. See Figure~\ref{fig:unbdd-face}.

Property~\ref{prop:regions} can be extended to the case of unbounded faces, as reflected in Figure~\ref{fig:unbdd-face}. Also compare with Figure~\ref{fig:Labels}.

\begin{figure}[h!]
	\centering
	\includegraphics[scale=0.6, page =3]{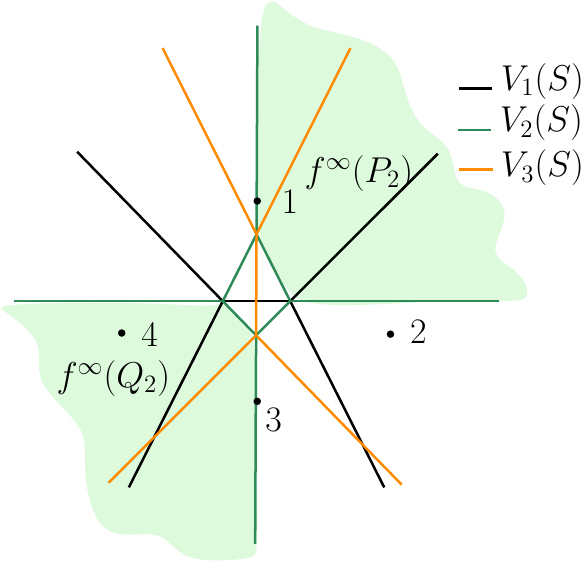}
	\caption{Inside $f^\infty(Q_2)$, there is an edge of $V_3(S)$ that crosses an edge of $V_1(S)$. But, inside $f^\infty(P_2)$, there is no edge of $V_3(S)$ that crosses an edge of $V_1(S)$. }
	\label{fig:cross-unbdd-face}
\end{figure}

The analogous statement of Property~\ref{prop:no-two-chains} for unbounded faces is not true; see Figure~\ref{fig:cross-unbdd-face}.

\section{Properties of the region $R_k(i)$}\label{sec:R_k}

Properties of geometric nature of the regions $R_k(i)$ are presented. Recall that $R_k(i)$ consists of all the points of the plane that have point $i \in S$ as one of their $k$ nearest neighbors from $S$, and $B_k(i)$ is the boundary of $R_k(i)$.

\begin{property}\label{prop:same-label}
    $R_k(i)$ forms a unique connected region. Its boundary $B_k(i)$ is formed by all the edges of $V_k(S)$ that have the label $i$, and this label is always inside $R_k(i)$. $B_k(i)$ is either a cycle, or one or more paths whose first and last edge are unbounded edges of $V_k(S)$. 
\end{property}
\begin{proof}
    Recall that, by Property~\ref{property:adjacentvoronoi}, proximity points are exchanged at the edges of $V_k(S)$. That is, two adjacent cells $f(P_k)$ and $f(Q_k)$ sharing the edge $\overline{b_{ij}}$ have $k-1$ common nearest neighbors from $S$ and they differ in $i\in P_k$ and $j\in Q_k$. The edge has the labels $i$ and $j$; by Property~\ref{property:cell-edge}, the label $i$ is inside $f(P_k)$ and the label $j$ is inside $f(Q_k)$. Therefore, the only edges that have the label $i$ are the edges incident to a cell in $R_k(i)$ and incident to a cell not in $R_k(i)$. Thus, the boundary of $R_k(i)$ is formed by all the edges of $V_k(S)$ which have the label $i$. Furthermore, by Property~\ref{property:cell-edge}, for all the edges of the boundary of $R_k(i)$, both the point $i$ and the region must be on the same side of the edge. This implies that if $R_k(i)$ was formed by more than one connected component, $i$ should be inside all of them, which is impossible. Hence $R_k(i)$ is connected.
    
    Note that each vertex of $V_k(S)$ is incident to either zero or two edges with the same label. Thus $B_k(i)$ is either a cycle, or one or more paths whose first and last edge are unbounded edges of $V_k(S)$.

\end{proof}

Edelsbrunner and Iglesias-Ham~\cite{EI18} proved that $R_k(i)$ is star-shaped. A slightly stronger statement holds, namely, we prove that the region $R_1(i)$ is contained in the kernel of $R_k(i)$ for every $k$. A polygon $P$ is star-shaped if there is at least one point $p\in P$ such that for any other point $q\in P$, the segment $\overline{pq}$ is strictly contained in $P$. The kernel of a star-shaped polygon $P$ is the set of points $p\in P$ such that $\overline{pq}$ is contained in $P$ for every $q\in P$. Equivalently, it is the intersection of all the half-planes defined by edges of $P$ which contain $P$.

\begin{property}\label{prop:star-shaped}
    For every $i\in S$, the region $R_k(i)$ of $V_k(S)$ is star-shaped, and the face $R_1(i)$ of $V_1(S)$ is contained in its kernel.
\end{property}
\begin{proof}
    Let $q$ be any point in $R_1(i)$. Then $q$ and $i$ lie to the same side of any bisector $b_{ij}$ for any $j\in S\setminus\left\lbrace i\right\rbrace$ because $i$ is the nearest neighbor of $q$ in $S$. Then $q$ is in the intersection of all the half-planes delimited by $b_{ij}$
    that contain $i$. This implies that $R_k(i)$ is star-shaped and that $q$ is in its kernel.
 \end{proof}

A vertex $v$ of $B_k(i)$ is reflex if the angle spanned by the two edges of $B_k(i)$ incident to $v$, inside $R_k(i)$, is greater than $\pi$. It is a convex vertex if this angle is less than $\pi$.

\begin{property}\label{prop:cycle-same-label}
The vertices of $B_k(i)$ that are incident to an edge of $V_k(S)$ lying in the interior (exterior) of $R_k(i)$ are of type II (type I) in $V_k(S)$. Moreover, for $k>1$, if $B_k(i)$ is a cycle, then it encloses at least three faces of $V_k(S)$. If $B_k(i)$ has $r$ reflex vertices, then it encloses at least $r$ faces of $V_k(S)$. 
\end{property}
\begin{proof}
    By Property~\ref{prop:same-label}, all the edges of $B_k(i)$ have the label $i$ inside $R_k(i)$; see Figure~\ref{fig:R32}. Hence, for a vertex $v\in B_k(i)$ that has an incident edge lying in the exterior of $R_k(i)$, there are two consecutive $i$ labels around $v$. Thus, the cyclic order of the labels of the edges around $v$ is $i, i, j, j, \ell, \ell$. Then, by Property~\ref{prop:cyclic-order}, $v$ must be of type I. For the remaining vertices, they have an incident edge lying in the interior of $R_k(i)$. This edge separates the two labels $i$. Thus, the cyclic order of the labels around $v$ must be $i, j, \ell, i, j, \ell$ and therefore, by Property~\ref{prop:cyclic-order}, $v$ must be of type II.

	\begin{figure}[h!]	
		\begin{center}
			\includegraphics[width=0.3\textwidth]{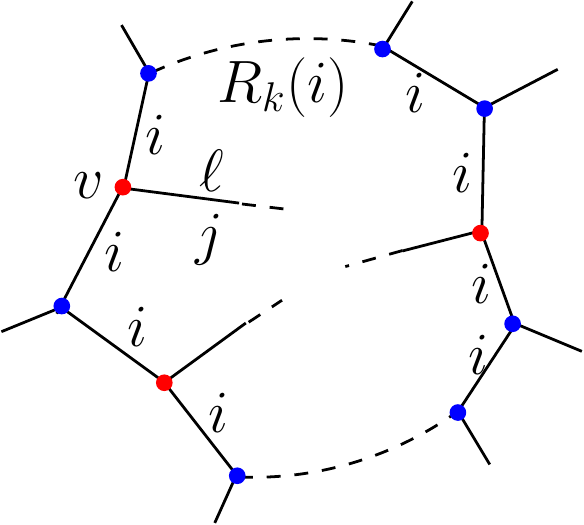}
			\caption{The vertices of $B_k(i)$ that have an incident edge lying in the interior of $B_k(i)$ are of type II. The remaining vertices of $B_k(i)$ are of type I.}		
			\label{fig:Flabel-in}
		\end{center}
	\end{figure}
	
	Let $B_k(i)$ be a cycle and $k\geq 2$. By Property~\ref{prop:no-two-chains}, it cannot happen that all vertices of a bounded cell of 
	$V_k(S)$ are of the same type or that the vertices of the same type are all consecutive in the boundary of the cell. Therefore, there are at least two edges of $V_k(S)$ in the interior of $B_k(i)$, which implies that $B_k(i)$ cannot be the boundary of a unique cell or of the union of two cells; see Figure~\ref{fig:Flabel-in}.
	
Assume then, that $B_k(i)$ has $r$ reflex vertices. Each reflex vertex $v$ of $B_k(i)$ has an incident edge 
$e_v$ that lies in the interior of $B_k(i)$, because otherwise the face incident to $v$ in $R_k(i)$ would not be convex. Note that each reflex vertex of $B_k(i)$ is of type II. When walking along the boundary $B_k(i)$ of $R_k(i)$ in clockwise order, we can assign to each reflex vertex of $B_k(i)$ the face of $R_k(i)$ that lies to the left of $e_v$, when orienting this edge towards the interior of $R_k(i).$ Then, to each reflex vertex we assign a different face in $R_k(i).$ It follows that $B_k(i)$	encloses at least $r$ cells of $V_k(S)$. 
	
\end{proof}

The next property follows from Property~\ref{prop:cycle-same-label}.

\begin{property}\label{prop:convex-reflex-vertex}
Let $v$ be a vertex in the boundary $B_k(i)$ of $R_k(i)$. If $v$ is of type I in $V_k(S)$, then it is a convex vertex of $B_{k}(i)$. If $v$ is of type II in $V_k(S)$, it is a reflex vertex of $B_{k}(i)$.
\end{property}

\begin{proof}
If $v$ is of type II, then it has an incident edge that lies in the interior of $R_k(i).$ Then, since all faces of $R_k(i)$ are convex, $v$ must be a reflex vertex. If $v$ is of type I, then it has an incident edge that lies in the exterior of $R_k(i).$ Using again that all faces of $R_k(i)$ are convex, $v$ must be a convex vertex. 
\end{proof}

\begin{property}\label{FRconnected}
    Let $f^b(P_k)$ be a bounded face of $R_k(i)$. Then $R_k(i)\setminus f^b(P_k)$ is connected.
\end{property}
\begin{proof}
For the sake of a contradiction, suppose that $R_k(i)\setminus f^b(P_k)$ is not connected. Let $F$ be a set of faces of $R_k(i)$ that forms one of the connected components of $R_k(i)\setminus f^b(P_k)$; see Figure~\ref{fig:ultima}.
 Necessarily, $F$ and $f^b(P_k)$ share a chain of edges $e_1, \dots, e_j $ and all the other edges of the boundary of $F$ are boundary edges of $R_k(i)$ having label $i$ inside $F$. Let $v_1$ and $v_j$ be the endpoints of the path $e_1, \dots, e_j $. In particular, $v_1$ and $v_j$ are boundary vertices of $R_k(i)$ and therefore incident to two edges with label $i$. For both vertices, one of the edges is a boundary edge of $F$. The edges in the chain $e_1, \dots, e_j $ do not have label $i$ because they are interior edges of $R_k(i)$. Then the edge of $f^b(P_k)$ incident to $v_1$, respectively $v_j$, not shared by $F$ has label $i$.
 By Property~\ref{prop:consecutive_labels} all the other edges of $f^b(P_k)$ different from $e_1,\dots, e_j$ have label $i$.
Hence, all the edges of $f^b(P_k)$ are either edges of $B_k(i)$ or edges shared with the component $F$, which implies that $R_k(i)= F \cup f^b(P_k)$ and $R_k(i)\setminus f^b(P_k)$ is connected. 

\begin{figure}[ht]
    \centering
    \includegraphics[scale=0.65]{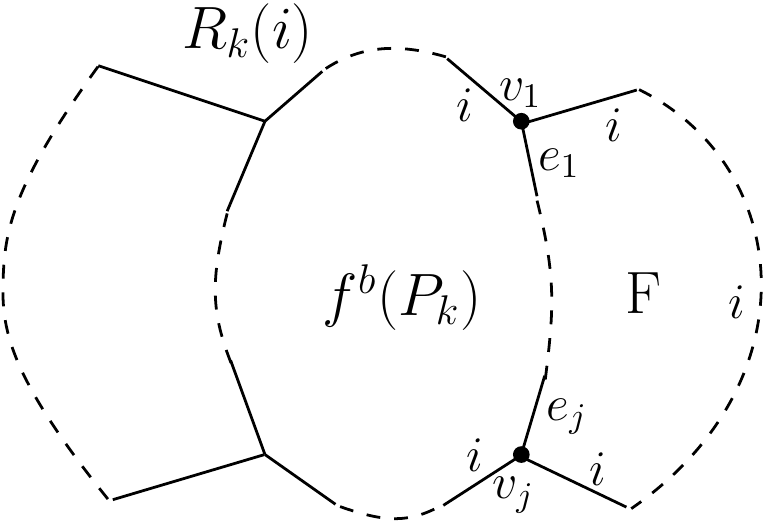}
    \caption{Illustration of the proof of Property~\ref{FRconnected}.}
    \label{fig:ultima}
\end{figure}
\end{proof}

\begin{property}\label{prop:nested}
For every $k>1$, and for every point $i \in S$, $R_{k-1}(i) \subset R_{k}(i).$ 
That is, the regions associated to the same point $i$ in Voronoi diagrams of consecutive orders are nested.
\end{property}
\begin{proof}

	\begin{figure}[h!]
	\centering
		\includegraphics[scale=0.55,page=2]{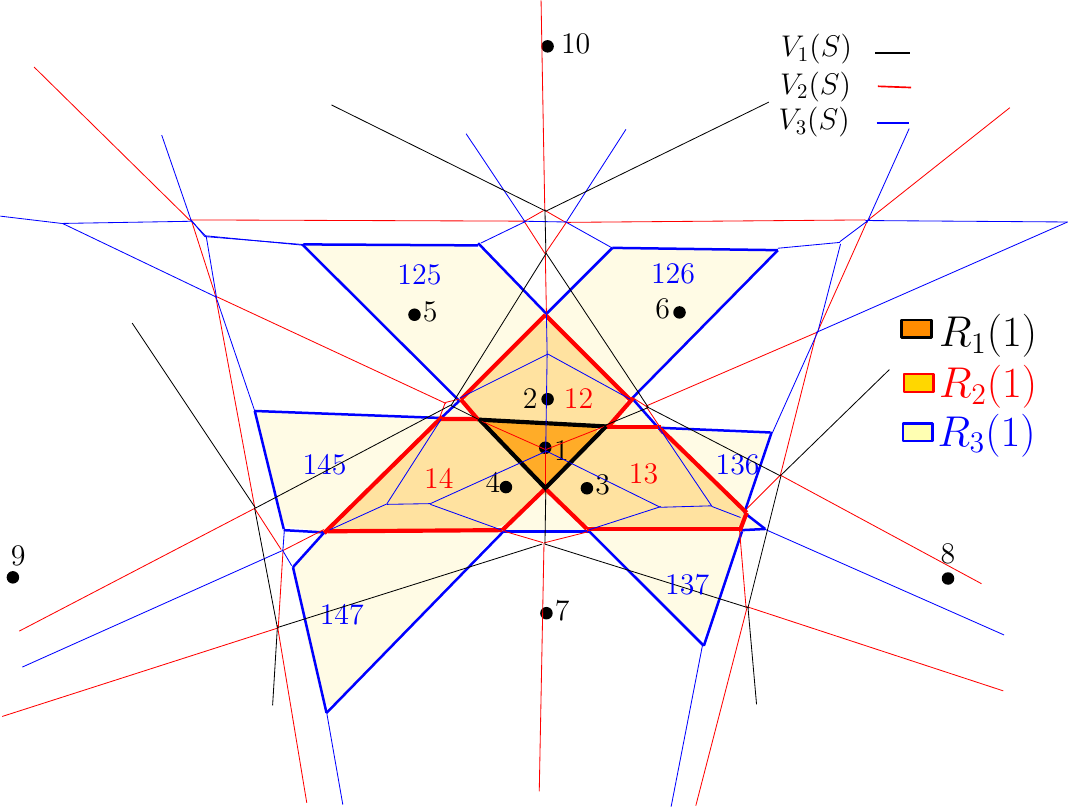}
	\caption{Nested regions $R_k(1)$ formed by the cells that have the point $1$ as one of the $k$ nearest neighbors, for $k\in\left\lbrace 1,2,3\right\rbrace$.}
	\label{fig:Fsame-label}
	
\end{figure}

Let $p\in R_{k-1}(i)$. Then $i$ is one of the $k-1$ nearest neighbors from $S$ of $p$. Hence, $i$ must be one of the $k$ nearest neighbors from $S$ of $p$. And therefore $p\in R_k(i)$. See Figure~\ref{fig:Fsame-label}.
\end{proof}

Let $\mathcal{A}(i)$ be the arrangement of the bisectors between point $i$ and another point of $S$.
Since $S$ is in general position, at most two bisectors intersect in a common point. 

Next, we show that $B_k(i)$ is determined by $B_{k-1}(i)$ using the arrangement $\mathcal{A}(i)$. This also provide a way to  construct the region $R_k(i)\backslash R_{k-1}(i)$ of a point $i \in S$ in the $k$-th nearest point Voronoi diagram of $S$ \cite{OBSC2000}.

\begin{property}
For $2\leq k\leq n-1$, $B_k(i)$
can be obtained from $B_{k-1}(i)$ by going through the arrangement of bisectors of $\mathcal{A}(i)$ in an orderly way. 
\end{property}
\begin{proof}
Any bisector $b_{ij}$ in $\mathcal{A}(i)$ is cut into $n-1$ open segments, delimited by the centers of the circles through $i$, $j$ and one of the $n-2$ remaining points of $S$. Note that two of the segments in $b_{ij}$ are unbounded, and that $\mathcal{A}(i)$ contains all the vertices and edges of $V_k(S)$ that contain the label $i$, for every $k$. In particular, $\mathcal{A}(i) = \bigcup_{k=1}^{n-1}B_k(i)$.  Let $\rm{v}$ be a vertex of $B_k(i)$. 
By Property~\ref{prop:cyclic-order}, $\rm{v}$ is incident to exactly four edges that have the label $i$. We describe the steps of an algorithm for computing $B_{k+1}(i)$ from $B_k(i)$. First, set $B_{k+1}(i)$ as an empty list of edges:
	\begin{enumerate}
		\item  For each convex vertex $\rm{v}$ of $B_k(i)$ we add  to $B_{k+1}(i)$ the two edges incident to $\rm{v}$ that are not in $B_k(i)$ and their endpoints. See Figure ~\ref{fig:construction}. If one of these vertices has degree $1$ in $B_{k+1}(i)$, we go to step $2$. Otherwise, if all the endpoints of these edges have degree $2$ in $B_{k+1}(i)$, we have obtained $B_{k+1}(i)$ and the algorithm stops.
		

		\item 
		Let $\overline{\rm{v}\rm{w}}$ be an edge of $B_{k+1}(i)$ such that $\rm{v}$ has degree $2$ in $B_{k+1}(i)$ and $\rm{w}$ has degree $1$. There are three possible edges of $\mathcal{A}(i)$ incident to $\rm{w}$ that we could add to $B_{k+1}(i)$. One of them is consecutive to $\overline{\rm{v}\rm{w}}$ along the same bisector and therefore cannot be an edge of $B_{k+1}(i)$. Among the two remaining edges, we add the one that is to the same side of the line containing $\overline{\rm{v}\rm{w}}$ as the point $i$. 	See Figure~\ref{fig:construction}. We repeat this step until there are no vertices of degree $1$ left. Then we have obtained $B_{k+1}(i)$.
		

\end{enumerate}	

We now justify the correctness of the algorithm.
By Property~\ref{prop:convex-reflex-vertex}, each convex vertex of $B_k(i)$ is a reflex vertex in $B_{k+1}(i)$, and they are the only reflex vertices in $B_{k+1}(i)$.
In step 1 we add these vertices and their incident edges to $B_{k+1}(i)$.
When the algorithm is finished, all the vertices of $B_{k+1}(i)$ must have degree $2$ since $B_{k+1}(i)$ is a unique cycle or one or more paths whose first and last edge are unbounded edges, see Property~\ref{prop:same-label}. But after step 1, there might still be vertices of degree $1$ because there are missing edges in $B_{k+1}(i)$. In step 2 we address those vertices.

It only remains to prove that, among the two possible edges that we can add to a vertex $\rm{w}$ of degree $1$ in step 2, we chose the edge that belongs to $B_{k+1}(i)$. But among the two edges only one is visible from point $i$, preserving the property that $R_k(i)$ is star-shaped, see Property ~\ref{prop:star-shaped}, and that is precisely the chosen edge in step 2.

 \begin{figure}[h!]
\centering
\subfloat[]{
\includegraphics[scale=0.43,page=1]{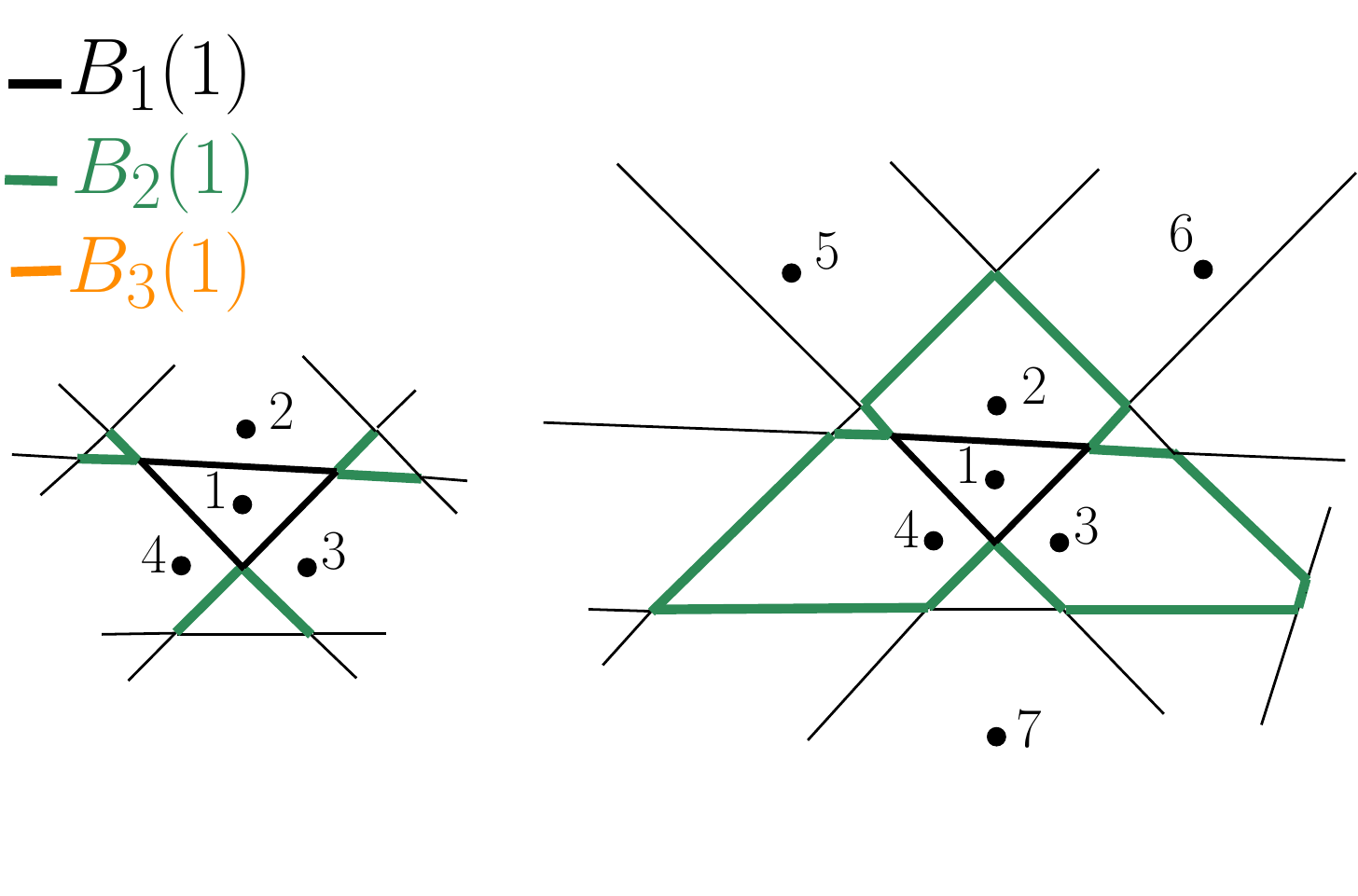}
}~
\subfloat[]{
\includegraphics[scale=0.43,page=2]{Falgorithm-Bki.pdf}
}
\caption{(a) Obtaining $B_2(1)$ (steps 1 and 2) from $B_1(1)$. (b)  Obtaining $B_3(1)$ (step 1) from $B_2(1)$.
}
\label{fig:construction}
\end{figure}

\end{proof}

\section{A small double cover of the edges of $V_k(S)$}\label{sec:double}

A cycle and path double cover of a graph $G$ is a collection of cycles and paths $\mathcal{C}$ such that every edge of $G$ belongs to precisely two elements of $\mathcal{C}$. 
We show in Property~\ref{prop:pathscycles} that $V_k(S)$ has a double cover $\mathcal{C}$ that is formed by the paths and cycles in $\cup_{i\in S}B_k(i)$, where $B_k(i)$ is the set of edges of $V_k(S)$ containing the same label $i$. 
The number of paths in such a double cover $\mathcal{C}$ is exactly the number of unbounded faces of $V_k(S)$, denoted by  $f_k^\infty$. Note that $f_k^\infty$ also is  the number of unbounded edges of $V_k(S)$. By next Property~\ref{lem:jedgeface}, $f_k^\infty$ also is equal to the number of $(k-1)$-edges $e_{k-1}$ of $S$. Proofs of this well-known observation can also be found in~\cite{E87,S98}.
It is known that  $e_{k-1}$ is at least $2k+1$~\cite{ELSS73} and at most $O(n\sqrt[3]{k})$~\cite{D98}.



\begin{property}\label{lem:jedgeface}
	For any set $S$ of $n$ points in general position,  $e_k=f_{k+1}^\infty.$
\end{property}
\begin{proof}
	Let $\overline{pq}$ be a $k$-edge of $S$ and consider a circle $C$ passing through $p$ and $q$ with center on the perpendicular bisector of $p$ and $q$, $b_{pq}$. When moving the center of $C$ along $b_{pq}$ towards infinity, $C$ approaches a half-plane defined by the supporting line of $\overline{pq}$, which divides $S\backslash\{p,q\}$ into two subsets of $k$ and $n-k-2$ points. Then, $b_{pq}$ contains an unbounded segment $s$ such that any circle $C$ passing through $p$ and $q$, with center on $s$, encloses $k$ points of $S$. Then, $s$ is an unbounded edge of $V_{k+1}(S)$. Each unbounded edge of $V_{k+1}(S)$ belongs to two unbounded faces of $V_{k+1}(S)$, and each unbounded face of $V_{k+1}(S)$ is incident to two unbounded segments. Then,  $e_k=f_{k+1}^\infty.$
\end{proof}


\begin{property}\label{prop:pathscycles}
$V_k(S)$,  which has $f_k^\infty$ unbounded faces, has an orientable cycle and path double cover consisting of $f_k^{\infty}$ paths and of at most $\max\{n-2k-1,2k-n-1\}$ cycles. The paths and cycles are given by the $B_k(i)$, for $1\leq i\leq  n-1.$ 
\end{property}

\begin{proof}
    For every $1 \leq i \leq n$, Property~\ref{prop:same-label} states that all the edges that have the label $i$ in $V_k(S)$ form a cycle or a set of paths. Each edge of $V_k(S)$ has two labels, and belongs therefore to exactly two cycles or paths. Since each vertex of $V_k(S)$ is incident to zero or two edges with label $i$, each path connects two unbounded edges of $V_k(S)$. Since each label $i$, corresponding to a point $i \in S$, is inside the corresponding region $R_k(i)$ (bounded by a cycle or by paths), we can orient all the edges of a cycle or path with label $i$ clockwise around point $i$, which is also contained in $R_k(i)$. This shows that the  cycle and path double cover is orientable. 
    It remains to show that the number of cycles is at most $\max\{n-2k-1,2k-n-1\}$. For $k\leq (n-2)/2$, we show that at least $2k+1$ different labels from $\{1,\dots,n\}$ appear as a label of some unbounded edge of $V_k(S)$, which implies that this double cover of $V_k(S)$ has at most $n-(2k+1)$ cycles. 
    An unbounded edge of $V_k(S)$ is an unbounded segment on some perpendicular bisector $b_{pq}$ of a segment $\overline{pq}$, with $p,q \in S$. This unbounded edge has labels $p$ and $q$. The segment $\overline{pq}$ is a $(k-1)$-edge of $S$; see Property~\ref{lem:jedgeface}. It is thus sufficient to show that at least $2k+3$ points of $S$ appear as points of $k$-edges of $S$; then at least $2k+3$ different labels appear in unbounded edges of $V_{k+1}(S)$. 
    Let $p$ be a point on the boundary of the convex hull of $S$. Sort radially the points of $S\backslash\{p\}$ around $p$. Let $p_1,\ldots, p_{n-1}$ be these ordered points. $\overline{pp_{k+1}}$ and $\overline{pp_{n-k-1}}$ are the two $k$-edges incident to $p$. Let then $q$ be any point from $\{p_1,\ldots,p_k,p_{n-k},\ldots,p_{n-1}\}.$ $pq$ is a $j$-edge for some $0\leq j <k,$ and it also is an $(n-j-2)$-edge. When rotating a line around $q$, the number of points on one side of the line changes from $j$ to $n-j-2$. The changes appear in steps of $1$, whenever the line passes through another point of $S$. Then $q$ is incident to a $k$-edge.
    Then, each point from the set of $2k+3$ points $\{p,p_1,\ldots,p_{k+1},p_{n-k-1},\ldots,p_{n-1}\}$ is incident to a $k$-edge. Consequently, the number of  cycles used in the double cover of $V_k(S)$ is at most $n-2k-1$. For $k=(n-1)/2$ or $k=n/2$ the same argument of rotating a line shows that each label from $\left\lbrace 1, \dots, n\right\rbrace$ appears on an unbounded segment of $V_k(S)$.
    For $k>n/2$, observe that $f_k^\infty=f_{n-k}^\infty$. More precisely, if an unbounded segment on a bisector $b_{pq}$ belongs to $V_k(S)$, then the other unbounded segment on $b_{pq}$ belongs to $V_{n-k}(S)$. Thus, the same labels appear in unbounded edges of $V_k(S)$ and in unbounded edges of $V_{n-k}(S)$. Then, at least $2(n-k)+1$ different labels appear in the unbounded edges of $V_k(S)$. Consequently, the number of cycles used in the double cover is at most $2k-n-1$.
\end{proof}

\begin{figure}[ht]
	\centering
		\includegraphics[scale=0.5]{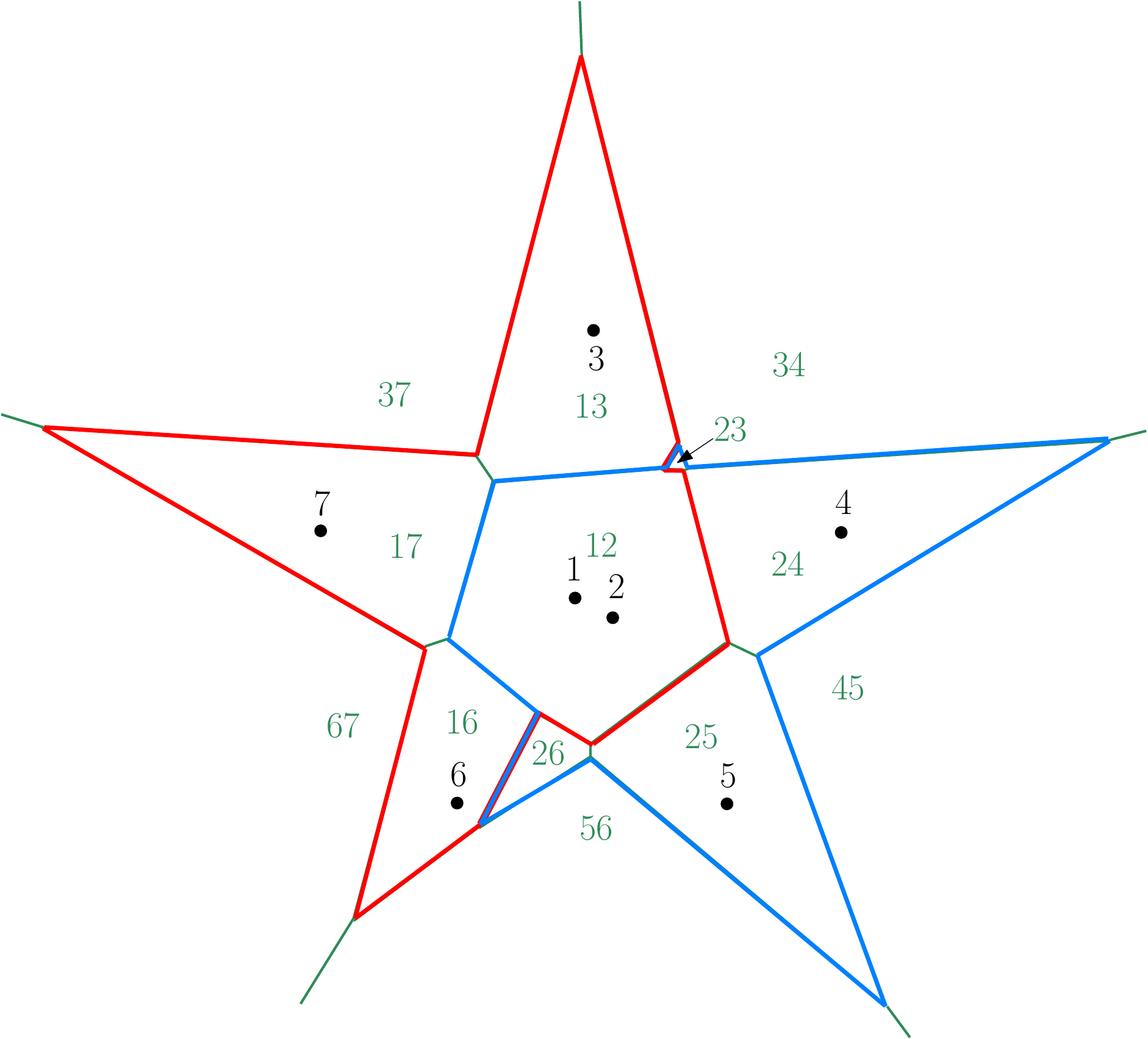}
	\caption{A set $S$ of $n=7$ points minimizing the number of unbounded regions of $V_2(S)$ among all point sets $S$ with $|S|\geq 5$. The boundaries of the regions $R_2(1)$ and $R_2(2)$ are cycles.}
	\label{cover}
	
\end{figure}

Note that this double cover is small compared to the number of vertices of $V_k(S)$, see~(\ref{number_vert}). The bound on the number of cycles in Property~\ref{prop:pathscycles} 
is attained for the point sets $S$ from~\cite{ELSS73} which have $2k+1$ $(k-1)$-edges. Such point sets consist of $2k+1$ points placed close to the vertices of a regular convex polygon and $n-2k-1$ points placed close to the center of the polygon. See Figure~\ref{cover}.
Property~\ref{prop:pathscycles} also shows that for $k=\lfloor{n/2}\rfloor$ and $k=\lceil{n/2}\rceil$, $V_k(S)$ has an orientable path double cover. This also holds for point sets in convex position and any value of $k$, see Property~\ref{prop:convexdouble}.
 
\begin{property}\label{prop:convexdouble}
    Let $S$ be a set of $n$ points in convex position. Then, $V_k(S)$ has an orientable path double cover consisting of $n$ paths.
\end{property}
\begin{proof}
    From the proof of Property~\ref{prop:pathscycles}, we know that $V_k(S)$ admits an orientable cycle and path double cover. We  show that every label from $\{1,\dots,n\}$ appears as a label of some unbounded edge, which implies that this double cover has no cycles. Indeed, every point of $S$ is incident to exactly two $(k-1)$-edges. By Property~\ref{lem:jedgeface}, each $(k-1)$-edge corresponds to an unbounded face of $V_k(S)$; more precisely, an unbounded segment on the bisector $b_{pq}$ perpendicular to a $(k-1)$-edge $\overline{pq}$ has labels $p$ and $q$. Then, each label appears in exactly two unbounded edges of $V_k(S)$. It follows that the double cover of $V_k(S)$ consists of exactly $n$ paths.
\end{proof}


\section{The number of vertices of $V_k(S)$}\label{sec:count_V_k(S)}

A formula for the number of vertices of $V_k(S)$ is well known~\cite{E87,L82,S98}. 
We present a new proof based on point movements and circles containing points in the following. 
Next Property~\ref{prop:vertices2} can also be found in~\cite{L03}, where it is proved using the formula for $V_k(S)$ from~\cite{E87,L82}. Recall that $e_j$ denotes the number of $j$-edges of $S$ and, by Property~\ref{lem:jedgeface}, equals the number of unbounded faces of $V_{j+1}(S)$. 

\begin{property}\label{prop:vertices2}
Let $c_k$ denote the number of circles through three points of $S$ that enclose exactly $k$ other points of $S$. Then
\begin{equation}\label{eq:ck}
c_k=(k+1)(2n-k-2)-\sum_{j=0}^{k} e_j
\end{equation}
\end{property}

\begin{proof}
The proof is based on continuous point moves, where at each moment only one point of the given set $S$ moves, and such that there never are more than four points cocircular and more than three points collinear. Each set of $n$ points in general position can be reached from any other set of $n$ points by performing such point moves. We show that Equation~(\ref{eq:ck}) remains valid throughout these movements. The value $c_k$ can only change when either (A) one point of $S$ leaves or enters a disk defined by three other points of $S$, or when (B) one point passes across the line spanned by two other points of $S$. In this second case, at the moment when the three points become collinear, the disk defined by the three points degenerates to a half-plane,  and then switches to a disk in the opposite half-plane; see Figure~\ref{fig:Fcircles}. We point out that for any $0 \leq j\leq k$, point moves across a line defined by two points $p_i$ and $p_j$ which form an $j$-edge of $S$  can modify the value of $c_k$. It was already observed in~\cite{R90} that the topological structure of $V_1(S)$ only changes at cases (A) and (B) when performing continuous point moves. The changes in $V_k(S)$ appear in the same way.\\
In case (A), observe that right before, and also right after, four points $P=\{p_a, p_b, p_c, p_d\}$ become cocircular, they form the vertices of a convex quadrilateral. We only need to examine how the number of points of $S$ inside a disk defined by three points from $P$ changes. This number can only change by $\pm 1$, depending on whether the fourth point of $P$ lies in the interior or in the exterior of disk defined by the other three points from $P$.  
The sign of the following determinant tells us, whether point $p_d$ is in the interior or in the exterior of the disk defined by $p_a, p_b,p_c$, see e.g.~\cite{GS85,R90}: 
\begin{equation}\label{eq:det}
\begin{vmatrix} x_{p_a} & y_{p_a} & x_{p_a}^2 + y_{p_a}^2 & 1\\
x_{p_b} & y_{p_b} & x_{p_b}^2 + y_{p_b}^2 & 1\\
x_{p_c} & y_{p_c} & x_{p_c}^2 + y_{p_c}^2 & 1\\
x_{p_d} & y_{p_d} & x_{p_d}^2 + y_{p_d}^2 & 1
\end{vmatrix}
\end{equation}
If $p_a,p_b,p_c$ defines a counter-clockwise oriented triangle, and if the determinant has positive sign, then $p_d$ is in the interior of the disk through $p_a,p_b,p_c.$ 
The determinant is zero iff the four points are cocircular. At this moment of cocircularity during the point movement, each of the four disks defined by $P$ contains the same number of points of $S$ in its interior. Using this determinant test, it is easy to verify that, before and after this movement, exactly two out of the four disks contain the fourth point of $P$ in the interior, also see~\cite{GS85}. Then, when performing a point move of case (A), the number $c_k$ does not change.\\ 
In case (B), a point $p_a$ of $S$ moves across the line spanned by points $p_b$ and $p_c$ of $S$. We can assume that $p_a$ crosses the interior of the segment $\overline{p_bp_c}$, as otherwise we could move one of $p_b$ and $p_c$ instead of $p_a$. Right before $p_a$ crosses $\overline{p_bp_c}$, $\overline{p_bp_c}$ is an $j$-edge of $S$ for some $0 \leq j\leq n-3$, and assume that before the movement $p_a$ is on the side of $\overline{p_bp_c}$ not containing these $j$ points. $\overline{p_bp_a}$ and $\overline{p_ap_c}$ are $(j+1)$-edges of $S$. See Figure~\ref{fig:Fcircles}. When moving $p_a$ across $\overline{p_bp_c}$, $\overline{p_bp_c}$ becomes an $(j+1)$-edge, and $\overline{p_bp_a}$ and $\overline{p_ap_c}$ become $j$-edges.
We get that $c_j$ decreases by one, $c_{n-j-3}$ increases by one, $e_j$ increases by one, $e_{j+1}$ decreases by one, $e_{n-j-2}$ increases by one, and $e_{n-j-3}$ decreases by one.
It is straightforward to verify that Equation~(\ref{eq:ck}) remains valid, independent of the value of $j.$\\
\begin{figure}[h!]
	\centering
	\subfloat[]{
		\includegraphics[scale=0.55,page=1]{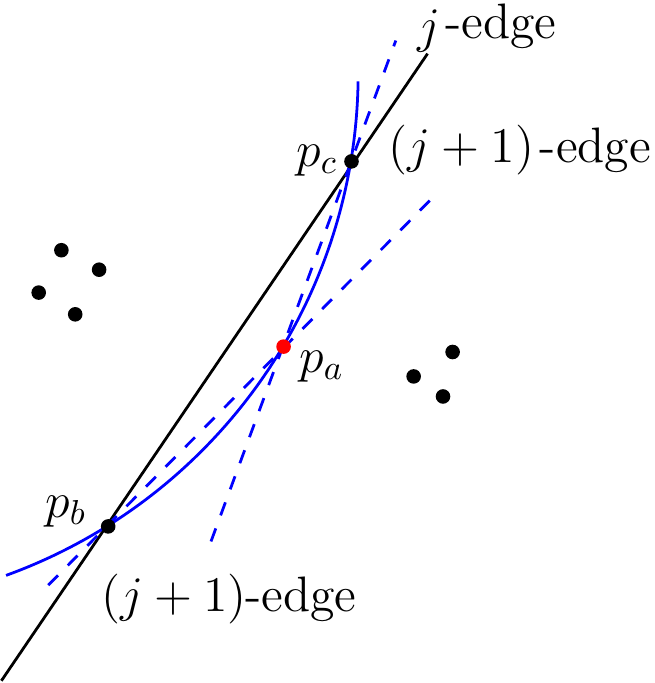}
	}~
	\subfloat[]{
		\includegraphics[scale=0.55,page=2]{Fcircles2.pdf}
	}~
	\subfloat[]{
		\includegraphics[scale=0.55,page=3]{Fcircles2.pdf}
	}
	\caption{When point $p_a$ is moved across the segment $\overline{p_bp_c}$, the disk defined by $p_a$, $p_b$ and $p_c$ in (a) degenerates to a half-plane in (b) and then switches to a disk in the opposite half-plane in (c). The change in the number of $j$-edges is also indicated.}
	\label{fig:Fcircles}
\end{figure}

It remains to show that Equation~(\ref{eq:ck}) is satisfied for some set of $n$ points. Let $S=\{p_1,\ldots, p_n\}$ with $p_i=(i,i^2)$, for $1 \leq i \leq n$.
Consider the disk defined by three given points $p_a,p_b, p_c$ of $S$, with $a<b<c$. To see if it contains another point $p_d$, we can use the determinant~(\ref{eq:det}), which in this case simplifies to 
$det=(c-d)(b-d)(a-d)(c-b)(c-a)(b-a)(a+b+c+d).$ Then, the sign of the determinant is positive only if $d<a$ or if $b<d<c.$ We now can count the number of disks defined by $S$ that contain $k$ points in their interior. These are all possibilities of choosing $p_a$ and $p_b$ such $a\leq k+1$ and $a+1\leq b \leq n-(k-a-1)$. The position of point $p_c$ is then determined. For each choice of $p_a$  there are $n-k-2$ choices for choosing point $p_b$; and we obtain that $c_k=(k+1)(n-k-2).$ Each point of $S$ is incident to two $j$-edges, for every $0\leq j \leq n-2.$ Then, $\sum_{j=0}^{k} e_j = (k+1)n$. Finally,
$c_k=(k+1)(2n-k-2)-\sum_{j=0}^{k} e_j$, as claimed. 
\end{proof} 

\begin{property}\label{prop:countvertices}
    $V_k(S)$ has $(4k-2)n-2k^2-\left(\sum_{j=0}^{k-1}e_j+ \sum_{j=0}^{k-2}e_j\right)$ vertices.
\end{property}
\begin{proof}
Let $v_k$ be the number of vertices of $V_k(S).$ By Property~\ref{property:points-types-edges}, the vertices of $V_{k}(S)$ are the centers of the circles through three points of $S$ that enclose exactly $k-1$ or $k-2$ points of $S$. Therefore
\begin{equation}
v_k=c_{k-1}+c_{k-2}.
\end{equation}
Also see~\cite{L03}.
Then, by Property~\ref{prop:vertices2}, $$v_k= (4k-2)n-2k^2-\left(\sum_{j=0}^{k-1}e_j+ \sum_{j=0}^{k-2}e_j\right).$$
\end{proof}

Since $\sum_{j=0}^{k}e_j\leq (k+1)n$, for $k<(n-2)/2$, see~\cite{AG86},  we get as a  consequence of Property~\ref{prop:countvertices}, that $V_k(S)$ has at least
\begin{equation}\label{number_vert}
  (2k-1)n-2k^2
\end{equation}
 vertices, for $k<(n-2)/2$. When $S$ is in convex
position, then $V_k(S)$ has this number of vertices; observe that for point sets in convex position we have $e_j= n$ for every $0\leq j\leq n-2$.\\

Another straight-forward consequence of Property~\ref{prop:countvertices} is the well-known formula for the number of faces of $V_k(S)$, see~\cite{L82}.
\begin{property}\label{prop:countfaces}
    $V_k(S)$ has $(2k-1)n-k^2+1- \sum_{j=0}^{k-2}e_j$ faces.
\end{property}
\begin{proof}
This property follows from Euler's relation among the number of vertices $|V|$, edges $|E|$, and faces $|F|$ of $V_k(S)$, see e.g.~\cite{D05, L03}:
\begin{equation}\label{eq:Euler}
|V|-|E|+|F|=1.
\end{equation}
Since $S$ is in general position, each vertex of $V_k(S)$ has degree $3$. By Property~\ref{lem:jedgeface}, the number of unbounded edges of $V_k(S)$ is $e_{k-1}.$ Then, $3|V|=2|E|-e_{k-1}.$ Then, Equation~(\ref{eq:Euler}) can be written as $|F|=\frac{|V|}{2}+\frac{e_{k-1}}{2}+1.$ It remains to substitute the formula for $|V|$ from Property~\ref{prop:countvertices}.
\end{proof}

\begin{property}\label{prop:difftypes}
The difference between vertices of type I and of type II of $V_k(S)$ is $2(n-k)-e_{k-1}.$

\end{property}
\begin{proof}
By Property~\ref{property:points-types-edges}, the number of vertices of type I in $V_k(S)$ is $c_{k-1}$, and  the number of vertices of type II in $V_k(S)$ is $c_{k-2}$. By Property~\ref{prop:vertices2}, $c_{k-1}-c_{k-2}=2(n-k)-e_{k-1}$.
\end{proof}


\section{Alternating hexagons}\label{sec:hexagons}
\begin{figure}[h!]
	\centering
		\includegraphics[scale=0.7]{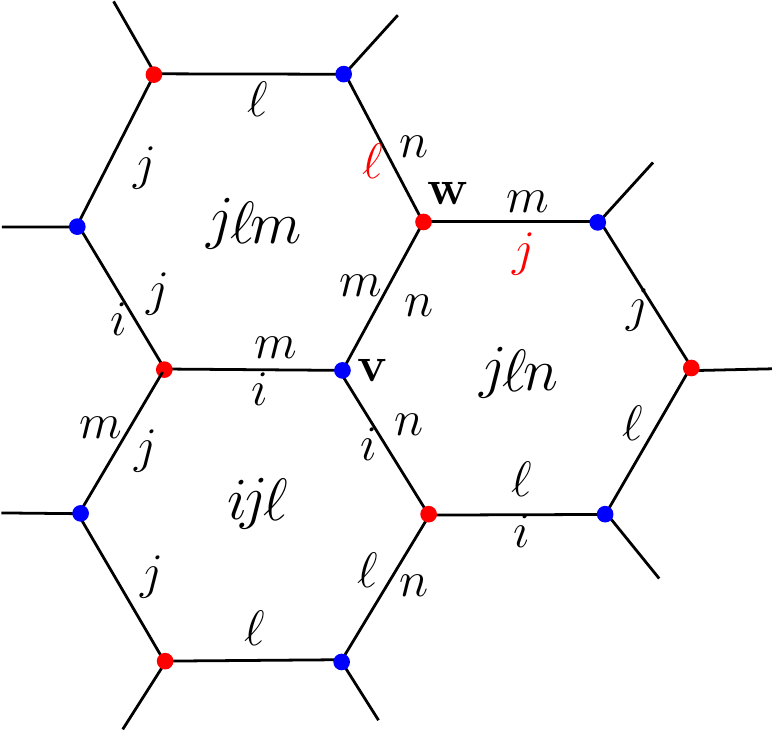}
	\caption{Three alternating hexagons incident to vertex $v$ of type I. The labels shown in red, corresponding to edges incident to $\rm{w}$, do not comply with Property~\ref{prop:cyclic-order}.}
	\label{fig:HEXAGONS}
\end{figure}

This section is on {\it{alternating hexagons}} and on certain configurations that cannot appear in $V_k(S)$ for small values of $k$.
A hexagonal face in $V_k(S)$ is alternating if its vertices alternate between type I and type II, see Figure~\ref{fig:HEXAGONS}. One might expect that a \lq\lq typical" face of $V_k(S)$ could be an alternating hexagon, because the average number of vertices in a face in a higher order Voronoi diagram is very close to six (Properties~\ref{prop:countvertices} and~\ref{prop:countfaces}), and because the numbers of vertices of type I and of type II in $V_k(S)$ do not differ a lot (Property~\ref{prop:difftypes}).  For $k=2$,  it follows from  
Property~\ref{prop:trianglefree} that $V_2(S)$ contains no alternating hexagons. We consider then $V_3(S)$. 

\begin{property}\label{prop:hexagons-around- typeI}
   Let $v$ be a vertex of type I in $V_3(S)$. Then, at most two of the three incident faces to $v$ are alternating hexagons.
\end{property}
\begin{proof}
    Let $\rm{v}$ be a vertex of type I in $V_3(S)$ and suppose that the three faces incident to $\rm{v}$ are alternating hexagons; see Figure~\ref{fig:HEXAGONS}. By Property~\ref{property:adjacentvoronoi} we can assume that these hexagons are defined by $P_3=\{i, j,\ell\}$, $P_3'=\{j, \ell, m\},$ and $P_3''=\{j,\ell,n\}.$ 
    By Property~\ref{prop:labels}, the labeling of one of these hexagons, say $f(P_3)$, determines the labeling of the other two hexagons. But then their common vertex $\rm{w}$ of type II does not comply with Property~\ref{prop:cyclic-order}. See  Figure~\ref{fig:HEXAGONS}. Therefore, there cannot be three alternating hexagons with a common vertex of type I in $V_3(S)$. 
\end{proof}
\begin{property}\label{prop:hexagons-V3}
    Any alternating hexagon in $V_3(S)$ is adjacent to at most three other alternating hexagons.  
\end{property}

\begin{proof}
    This follows directly from Property~\ref{prop:hexagons-around- typeI}.
\end{proof}

\begin{figure}[h!]
    \centering
    \includegraphics[scale = 0.45]{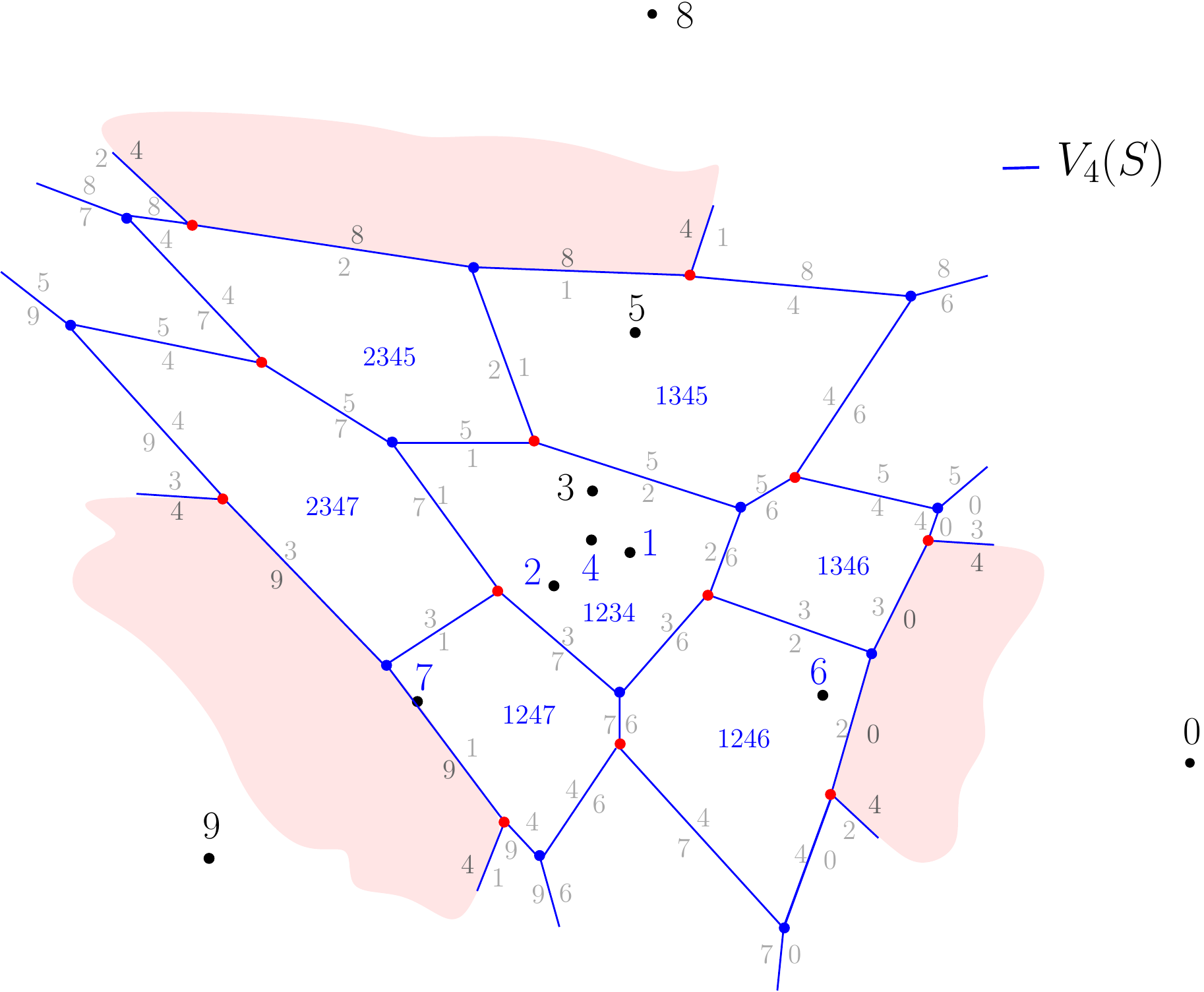}
    \caption{Example of a set $S = \left\lbrace 0, 1,\dots, 9\right\rbrace$ 
    such that in $V_4(S)$, all the faces adjacent to the alternating hexagon $f^b(\left\lbrace 1,2,3,4\right\rbrace)$ are also alternating hexagons. The light red faces cannot be alternating hexagons.}
    \label{fig:4alter}
\end{figure}

For $k=4$, in $V_4(S)$ it is possible that all the adjacent faces to an alternating hexagon are alternating hexagons, see Figure~\ref{fig:4alter}. The figure also shows that a second layer of alternating hexagons around the central face $f^b(\left\lbrace 1,2,3,4 \right\rbrace$ is not possible. The light red regions in Figure~\ref{fig:4alter} cannot be alternating hexagons because all the edges with the same label must be consecutive in a face. 

\begin{figure}[h!]
    \centering
    \includegraphics[scale=0.6]{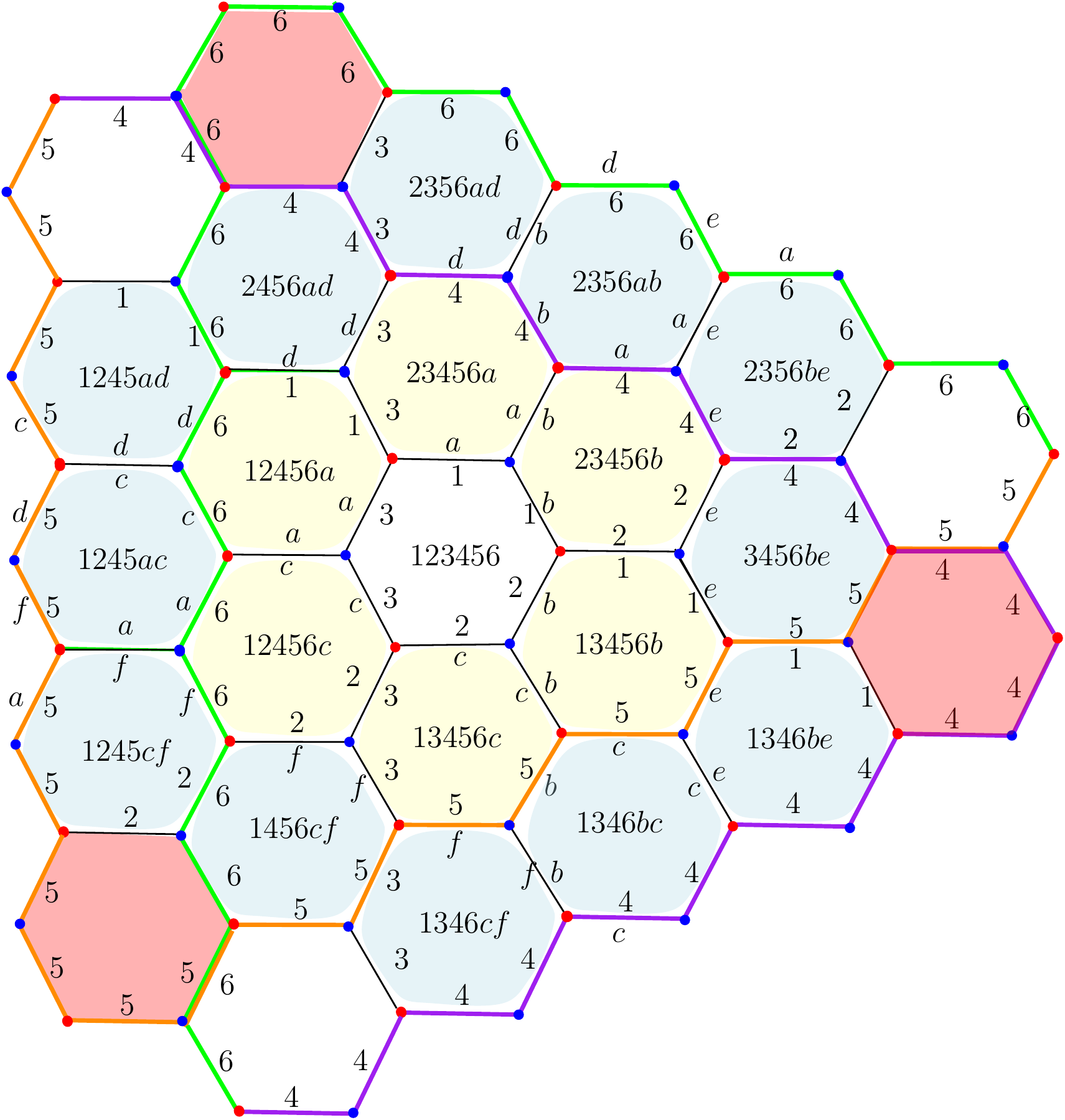}
    \caption{Labeling of alternating hexagons in $V_6(S)$. Two consecutive layers around the central face $f^b(\{1,2,3,4,5,6\})$ satisfy the labeling rules. The three red faces do not comply with the vertex rule.} 
    \label{fig:malla-hexa}
\end{figure}

For larger values $k$, the labeling rules can be satisfied for more layers of alternating hexagons surrounding a given alternating hexagon. In Figure~\ref{fig:malla-hexa} for $V_6(S)$ it is possible to form two consecutive layers (in yellow and blue) around the central alternating hexagon. However, it is not possible to form a third layer, independent of the labeling. See Figure~\ref{fig:malla-hexa} for one possible labeling. We leave for future investigation whether this labeling is realizable for some point set.

We finally present a property on alternating hexagons in $R_k(i)$.
\begin{property}
    For every $i\in S$, the faces forming the region $R_k(i)$ cannot all be alternating hexagons.
\end{property}
\begin{proof}
Note that if $R_k(i)$ is formed only by alternating hexagons, then there is always a vertex $v$ of type II in $B_k(i)$ without an incident edge lying in the interior of $R_k(i)$. By Property~\ref{prop:cyclic-order}, the two edges of $B_k(i)$ incident to $v$ cannot both have the label $i$ inside $R_k(i)$.
\end{proof}

\section{Acknowledgments}

M. Claverol and C. Huemer were supported by PID2019-104129GB-I00/ AEI/ 10.13039/501100011033 of the Spanish Ministry of Science and Innovation. M. Claverol was also supported by Gen. Cat. DGR 2017SGR1640 and C. Huemer was also supported by Gen. Cat. DGR 2017SGR1336.
A. Martínez-Moraian was funded by the predoctoral contract PRE2018-085668
of the Spanish Ministry of Science, Innovation, and Universities and was partially supported by project
PID2019-104129GB-I00/AEI/10.13039/501100011033. 

	\begin{minipage}[l]{0.3\columnwidth}
		\includegraphics[trim=10cm 6cm 10cm 5cm,clip,scale=0.15]{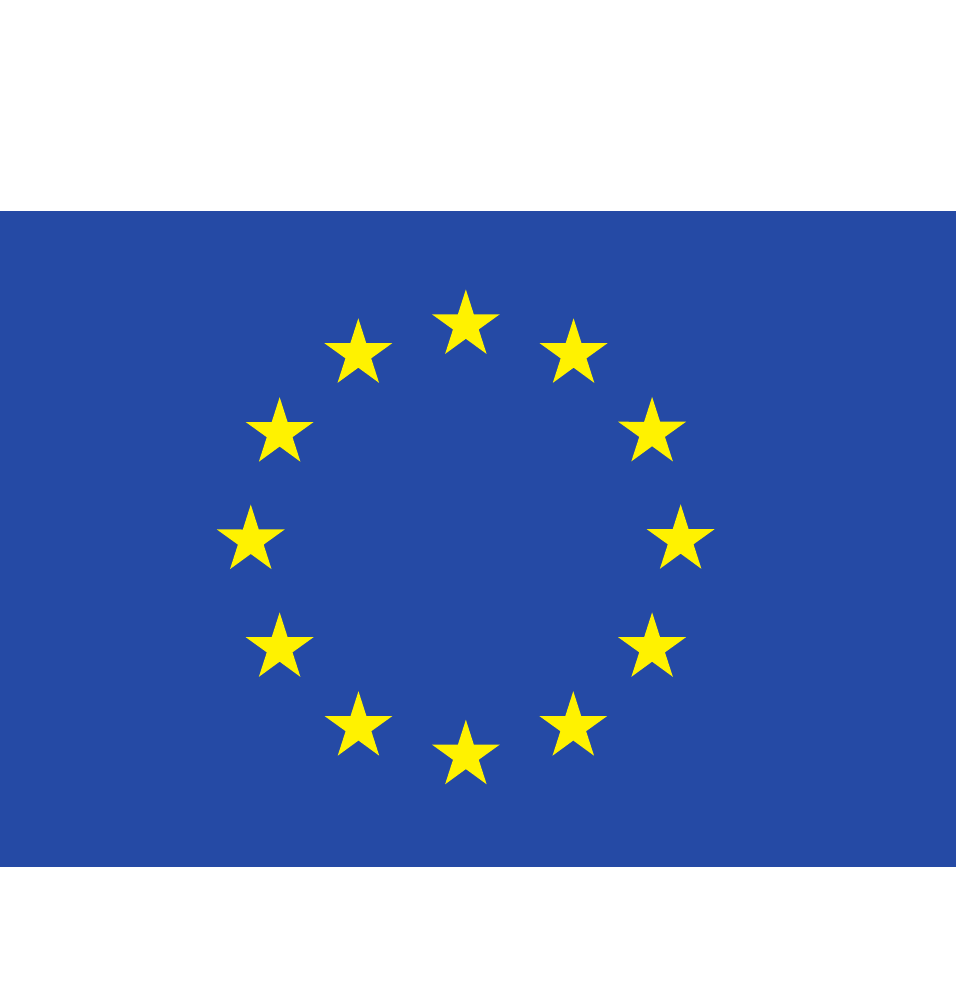}
	\end{minipage}
	\hspace{-1cm}
	\begin{minipage}[l][1cm]{0.75\columnwidth}
		This work has received funding from the European Union's Horizon 2020
		research and innovation programme under the Marie Sk\l{}odowska-Curie
		grant agreement No 734922.
	\end{minipage}

\newpage
{}

\end{document}